\documentclass[12pt,thmsa]{article}
\usepackage{amsmath}
\usepackage{graphicx}
\usepackage{amsfonts}
\usepackage{amssymb}
\usepackage{mathrsfs}

\newtheorem{theorem}{Theorem}

\newtheorem{lemma}{Lemma}

\newtheorem{remark}{Remark}

\begin{document}

\begin{center}
{\Large \textbf{An independence test for functional variables based on kernel normalized cross-covariance operator}}

\bigskip

Terence Kevin  MANFOUMBI DJONGUET and  Guy Martial  NKIET 

\bigskip

\textsuperscript{}URMI, Universit\'{e} des Sciences et Techniques de Masuku,  Franceville, Gabon.

\bigskip

E-mail adresses : tkmpro95@gmail.com,    guymartial.nkiet@univ-masuku.org.

\bigskip
\end{center}

\noindent\textbf{Abstract.}  We propose an independence test for random variables valued into metric spaces by using a test statistic obtained from appropriately  centering and rescaling the squared Hilbert-Schmidt norm of the usual empirical estimator of normalized cross-covariance operator.  We then get asymptotic normality of this statistic  under independence hypothesis, so leading to a  new  test for independence  of  functional random variables. A   simulation study that allows to compare   the proposed test   to  existing ones is provided.
\bigskip


\section{Introduction\label{intro}}

\noindent Testing the independence between two random variables has always been an important concern in statistics, and  many methods have been developed to do this in the case of real or multidimensional random variables. In the case where the concerned variables  are functional, one finds within Functional Data  Analysis (FDA) which is a  branch of statistics allowing   to set up methods for processing data being digitized points of   curves, some methods making it possible to test the independence or the lack  of  relationship  between the concerned variables. This is how  non-correlation tests between functional variables were proposed in \cite{aghoukeng10} and    \cite{kokoszka08}      whereas  independence tests were introduced in    \cite{gorecki20},   \cite{lai21}  and  \cite{manfoumbi22}. More specifically, G\'orecki et al.\cite{gorecki20}  proposed the use of the  Hilbert-Schmidt Independence Criterion (HSIC), introduced by Gretton et al. \cite{gretton05a},  for testing independence of   multivariate functional data, Lai et al. \cite{lai21}  introduced the angle covariance to characterize independence and proposed  a permutation based test whereas Manfoumbi Djonguet and Nkiet \cite{manfoumbi22}  proposed a modification of a  HSIC naive estimator that yields asymptotic normality under null hypothesis, so allowing a simple test.    HSIC is one of the dependence measures proposed within the framework of kernel-based methods, that is methods based on kernel embeddings of probability measures in a   reproducing kernel Hilbert space (RKHS). These are powerful methods allowing to work with high-dimensional and structured data (e.g., \cite{harchaoui13}), including functional data.  Other kernel-based dependence measures have been introduced in various works, notably in \cite{fukumizu07b} and \cite{gretton05b}. Particularly, Fukumizu et al. \cite{fukumizu07b} introduced a dependence measure based on the    so-called  normalized cross-covariance operator (NOCCO) defined in \cite{fukumizu07a}. More specifically, they used the  squared Hilbert-Schmidt norm of the  empirical estimator of  NOCCO for measuring independence and introduced a permutation test of independence based on the measure thus proposed. To the best of our knowledge,  there do not exist other studies of this dependence measure which is however very interesting since, as it is proved in \cite{fukumizu07b}, it characterizes independence when the related kernels are characteristic ones.  As a  permutation test has  the disadvantage of requiring high computational costs,   it is of interest to perform rather the test by using  the   asymptotic distribution of  a test statistic based on the  aforementioned measure under independence  hypothesis.  This is the approach we take    in this paper. More specifically, we introduce   a test statistic obtained from appropriately  centering and rescaling the squared Hilbert-Schmidt norm of the   empirical estimator of NOCCO, and we get   asymptotic normality  for  this statistic  under independence hypothesis. This allows to propose   a  new  test for independence  of   random variables valued into metric spaces, including functional random variables. The rest of the paper is organized as follows.      Section 2 is devoted to the introduction of basic notions that are used in the paper. In   Section 3, the NOCCO-based dependence measure introduced in \cite{fukumizu07b} is recalled, then we defined a test statistic based on it for which we obtain an asymptotic normality under the independence hypothesis and specified assumptions. We then show how this test statistic can be computed in practice, especially when functional data are available.   A simulation study on functional data that allows to compare   the proposed test   to those  of Gretton et al. \cite{gretton05a} and Fukumizu  et al. \cite{fukumizu07b}  is  given in Section 4. All the proofs are postponed in Section 5.
\section{Preliminary notions}
\noindent This section is devoted to recalling  the notion of reproducing kernel hilbert space (RKHS) and defining  some elements related to it that are useful in this paper. For more details on RKHS and its use in probability and statistics, one may refer to \cite{berlinet04}.
 Letting  $(\mathcal{X},\mathcal B_{\mathcal{X}})$ be a measurable space, where $(\mathcal{X},d)$ is a metric space  and $\mathcal B_{\mathcal{X}}$ is the corresponding Borel $\sigma$-field, we consider a Hilbert space $\mathcal{H}_\mathcal{X}$ of functions from $\mathcal{X}$ to $\mathbb{R}$, endowed with an inner product $<\cdot,\cdot>_{\mathcal{H}_\mathcal{X}}$. This space is said to be a RKHS   if there exists a kernel, that is a symmetric positive semi-definite function $K:\mathcal{X}^2\rightarrow \mathbb{R} $, such that for any $f\in \mathcal{H}_\mathcal{X}$ and any $x\in\mathcal{X}$, one has $K(x,\cdot )\in \mathcal{H}_\mathcal{X}$ and  $f(x)=<f,K(x,\cdot)>_{\mathcal{H}_\mathcal{X}}$. The feature map $\Phi\,:\,x\in\mathcal{X}\mapsto K(x,\cdot)\in \mathcal{H}_\mathcal{X}$ characterizes $K$ since $
K(x,y)=<\Phi(x),\Phi(y)>_{\mathcal{H}_\mathcal{X}}
$
for any $(x,y)\in\mathcal{X}^2$.    Let   $X$ be a random variable taking values in $\mathcal{X}$ and with probability distribution $\mathbb{P}_X$. If $\mathbb{E}(\Vert \Phi(X)\Vert_{\mathcal{H}_\mathcal{X}})<+\infty$, the mean element $m_X$ associated with  $X$ is defined for all functions $f\in \mathcal{H}_\mathcal{X}$ as the unique element in  $\mathcal{H}_\mathcal{X}$ satisfying
\begin{eqnarray*}
<m_X,f>_{\mathcal{H}_\mathcal{X}}=\mathbb{E}\left(f(X)\right)=\int_{\mathcal{X}}f(x)\,\,d\mathbb{P}_X(x).
\end{eqnarray*}
Furhermore, if  $\mathbb{E}\left(\Vert \Phi(X)\Vert_{\mathcal{H}_\mathcal{X}}^2\right)<+\infty$, we can define the covariance operator associated to $X$ as the unique operator $\Sigma_X$ from $\mathcal{H}_\mathcal{X}$ to itself such  that, for any pair $(f,g)\in\mathcal{H}_\mathcal{X}^2$, one has
\begin{eqnarray*}
<f,\Sigma_Xg>_{\mathcal{H}_\mathcal{X}}=\textrm{Cov}\left(f(X),g(X)\right)=\mathbb{E}\left(f(X)g(X)\right)-\mathbb{E}(f(X))\,\mathbb{E}(g(X)).
\end{eqnarray*}
It is very important to note that if   $\Vert K\Vert_\infty<+\infty$, where 
\[
\Vert K\Vert_\infty:=\sup_{(x,y)\in\mathcal{X}^2}K(x,y), 
\]
then   $m_X$  and   $\Sigma_X$ exist. They can also be expressed as
\[
m_X=\mathbb{E}\left(K(X,\cdot)\right)\,\,\,\textrm{ and }\,\,\,
\Sigma_X=\mathbb{E}\left(K(X,\cdot )\otimes K(X,\cdot)\right)-m_X\otimes m_X,
\]
where $\otimes$ is the tensor product such that, for any pair of vectors  $(x,y)$, the product $x\otimes y$ is the linear map defined by  $(x\otimes y)(h)=<x,h>\,y$ for all $h$. Now, let $Y$ be another random variable with values in a measurable space  $(\mathcal{Y},\mathcal B_{\mathcal{Y}})$, where $\mathcal{Y}$ is a metric space  and $\mathcal B_{\mathcal{Y}}$ is the corresponding Borel $\sigma$-field. We consider a RKHS $\mathcal{H}_\mathcal{Y}$ of functions from $\mathcal{Y}$ to $\mathbb{R}$, endowed with an inner product $<\cdot,\cdot>_{\mathcal{H}_\mathcal{Y}}$, and associated with a kernel  $L:\mathcal{Y}^2\rightarrow \mathbb{R} $. If  $\mathbb{E}\left(\Vert \Psi(X)\Vert_{\mathcal{H}_\mathcal{Y}}^2\right)<+\infty$, where $\Psi=L\left( Y,\cdot\right)$ and $\Vert\cdot\Vert_{\mathcal{H}_\mathcal{Y}}$ is the norm associated with $<\cdot,\cdot>_{\mathcal{H}_\mathcal{Y}}$, we can define the mean element $m_Y$ and  the covariance operator  $\Sigma_Y$ as indicated above.  We can also define the cross-covariance operator of $X$ and $Y$, denoted by $\Sigma_{XY}$, as   the unique linear operator   from $\mathcal{H}_{\mathcal{Y}}$ to $\mathcal{H}_{\mathcal{X}}$ satisfying
$$ \langle f , \Sigma_{XY}g \rangle_{\mathcal{H}_{\mathcal{X}}} = cov(f(X),g(Y))  $$
for all $(f,g) \in \mathcal{H}_{\mathcal{X}} \times \mathcal{H}_{\mathcal{Y}}$. It can also be expressed as
\[
\Sigma_{XY}=\mathbb{E}\bigg(\Big(L(Y,\cdot)-m_Y\Big)\otimes \left(K(X,\cdot)-m_X\right)\bigg)=\mathbb{E}\left(L(Y,\cdot )\otimes K(X,\cdot)\right)-m_Y\otimes m_X.
\]
It is well known (see, e.g., \cite{fukumizu07a,fukumizu07b})  that $\Sigma_{XY}$ is  an Hilbert Schmidt operator and that there exists a unique bounded  operator $V_{XY}\,:\,\mathcal{H}_\mathcal{Y}\rightarrow \mathcal{H}_\mathcal{X}$, called the normalized cross-covariance operator (NOCCO), satisfying
\begin{equation}\label{deco}
\Sigma_{XY} = \Sigma_X ^{1/2} V_{XY} \Sigma_Y ^{1/2}. 
\end{equation}
From now on, we assume that $V_{XY}$ is a Hilbert-Schmidt operator. Note that, Fukumizu et al. \cite{fukumizu07a} discussed conditions that ensure  this later property and  they show in Theorem 3 that  a sufficient condition is the fact that the mean square contingency is finite, what is very natural when two different random variables are considered.  As we will see later,  $V_{XY}$  measures  the dependence between $X$ and $Y$ and can, therefore, be the basis of a procedure for testing for independence.   
\section{NOCCO and independence testing}
\noindent In this section, we consider a dependence measure based on NOCCO that permits to consider a test of independence, that is the test of the hypothesis $\mathscr{H}_0\,:\,X\perp \!\!\! \perp Y$, where $\perp \!\!\! \perp$ denotes stochastic  independence, against the alternative hypothesis $\mathscr{H}_1$ meaning that $X$ and $Y$ are not independent. Then, a consistent estimator of this measure is proposed as test statistic, and we get its asymptotic normality under $\mathscr{H}_0$.
\subsection{A NOCCO-based dependence measure}
 \noindent One of the most famous RKHS-based dependence  measures is the Hilbert-Schmidt independence criterion (HSIC) introduced in \cite{gretton05a}. It  is equal to  $\Vert \Sigma_{XY}\Vert^2_{\textrm{HS}}$, where   $\Vert \cdot\Vert_{\textrm{HS}}$  denotes the Hilbert-Schmidt norm, and is known to fully characterize independence property. Indeed, when $K$ and $L$ are characteristic kernels $\mathscr{H}_0$ is equivalent to  $ \Sigma_{XY} =0$. That is why estimate  of HSIC was used for testing for independence (\cite{gretton05a,gretton08}).  Unfortunely, asymptotic normality under $\mathscr{H}_0$ could not be obtained for  the used test statistic  but rather convergence in distribution  to an infinite sum of chi-squared distributions, what is inusable for performing the test and   requires the use of methods such as permutation methods which have high computation costs. Later, Fukumizu et al. \cite{fukumizu07b}   proposed  to use instead the dependence measure $N(X,Y)=\Vert V_{XY}\Vert^2_{\textrm{HS}}$; it also characterizes independence since from \eqref{deco} it is clear that,  when $K$ and $L$ are characteristic kernels,  $\mathscr{H}_0$ is equivalent to  $N(X,Y)=0$. In order to achieve an independence test, we have to consider a consistent estimator of $N(X,Y)$.
\subsection{Test statistic and asymptotic normality}
\noindent Let $\{(X_i,Y_i)\}_{1\leq i\leq n}$ be a i.i.d. sample of $(X,Y)$;  the   empirical counterparts of $m_X$ and $m_Y$ are respectively given by: 
\begin{equation*}\label{empirm}
\overline{K}_n=\frac{1}{n} \sum_{i=1}^{n}K(X_i,\cdot)\,\,\,\textrm{ and }\,\,\,
\overline{L}_n=\frac{1}{n} \sum_{i=1}^{n}L(Y_i,\cdot).
\end{equation*}
We then consider the following estimators  of $\Sigma_{X}$, $\Sigma_{Y}$ and  $\Sigma_{XY}$ defined as
\begin{equation*}\label{empirv1}
\widehat{\Sigma}_X=\frac{1}{n} \sum_{i=1}^{n}\left(K(X_i,\cdot)-\overline{K}_n\right)\otimes\left(K(X_i,\cdot)-\overline{K}_n\right)
=\frac{1}{n} \sum_{i=1}^{n} K(X_i,\cdot)\otimes K(X_i,\cdot) -\overline{K}_n\otimes\overline{K}_n,
\end{equation*}
\begin{equation*}\label{empirv1}
\widehat{\Sigma}_Y=\frac{1}{n} \sum_{i=1}^{n}\left(L(Y_i,\cdot)-\overline{L}_n\right)\otimes\left(L(Y_i,\cdot)-\overline{L}_n\right)
=\frac{1}{n} \sum_{i=1}^{n} L(Y_i,\cdot)\otimes L(Y_i,\cdot) -\overline{L}_n\otimes\overline{L}_n,
\end{equation*}
and
\begin{equation*}\label{empirv1}
\widehat{\Sigma}_{XY}=\frac{1}{n} \sum_{i=1}^{n}\left(L(Y_i,\cdot)-\overline{L}_n\right)\otimes\left(K(X_i,\cdot)-\overline{K}_n\right)
=\frac{1}{n} \sum_{i=1}^{n} L(Y_i,\cdot)\otimes K(X_i,\cdot) -\overline{L}_n\otimes\overline{K}_n.
\end{equation*}
From them, we want to define an estimator of $V_{XY}$. A natural way to do this would be to invert  $\widehat{\Sigma}_X^{1/2}$ and  $\widehat{\Sigma}_Y^{1/2}$ but these operators are compact ones and, therefore,   do  not have inverses. So, as in \cite{fukumizu07b}, we use a regularization approach:  we consider a sequence $\{\gamma_n\}_{n\geq 1}$  of strictly positive numbers such that $\lim_{n\rightarrow +\infty}(\gamma_n)=  0$, and we estimate  $V_{XY}$ by the random operator
$$ \widehat{V}_{XY}    = (\widehat{\Sigma}_{X} + \gamma_n\mathbb{I}_\mathcal{X})^{-1/2}\,\widehat{\Sigma}_{XY}   (\widehat{\Sigma}_{Y} +\gamma_n\mathbb{I}_\mathcal{Y})^{-1/2}, $$
where  $\mathbb{I}_\mathcal{X}$ (resp.  $\mathbb{I}_\mathcal{Y}$) denotes the identity operator of $\mathcal{X}$ (resp. $\mathcal{Y}$). Then, we take as estimator of $N(X,Y)$ the statistic
$$ \widehat{N}_n (X,Y) = \Vert\widehat{V}_{XY}  \Vert_{\textrm{HS}} ^2. $$
 Theorem 5  in \cite{fukumizu07b}  establishes consistency of this estimator   when    $\gamma_n ^{-3} n^{-1} \rightarrow 0$ as $n\rightarrow +\infty$. So,    $ \widehat{N}_n (X,Y)$ can be used as test statistic for the problem of testing for $\mathscr{H}_0$ against  $\mathscr{H}_1$. It is, therefore, necessary to get asymptotic distribution under $\mathscr{H}_0$ of this statistic. In fact, we will obtain asymptotic normality after recentering and rescalling it. As in  \cite{harchaoui08}, for a given compact operator $A$ with decreasing eigenvalues $\lambda_p(A)$, $p\geq 1$, we define the quantity $d_r(A,\gamma)$ for all $r\in\{1,2\}$ and $\gamma\in]0,+\infty[$  as
\begin{equation*}\label{dr}
d_r(A,\gamma)=\bigg\{\sum_{p=1}^{+\infty}\left(\lambda_p(A)+\gamma\right)^{-r}\lambda_p(A)^r\bigg\}^{1/r}.
\end{equation*}
Then, for $r\in\{1,2\}$, we put
\begin{equation*}\label{drchap}
\widehat{D}_{r,n}=d_r(\widehat{\Sigma}_X,\gamma_n)\,d_r(\widehat{\Sigma}_Y,\gamma_n)
\end{equation*}
and we consider the  centered and rescaled version of  $ \widehat{N}_n (X,Y)$ given by
\begin{equation*}
\widehat{T}_n=\frac{n\widehat{N}_n (X,Y)-\widehat{D}_{1,n}}{\sqrt{2}\,\widehat{D}_{2,n}}.
\end{equation*}
In order to get the asymptotic distribution of  $\widehat{T}_n$, we make the following assumptions:

\bigskip

\noindent $(\mathscr{A}_1)$:  the kernels $K$ and  $L$ are universal and satisfy $\Vert K\Vert_{\infty}   <+ \infty $ and $\Vert L\Vert_{\infty}   <+ \infty $;

\bigskip

\noindent $(\mathscr{A}_2):$ the sequence of decreasing eigenvalues  $\{\lambda_p\}_{p\geq 1}$ of $\Sigma_X$ is such that  $\sum_{p=1}^{+\infty}\lambda_p^{1/2} <+\infty$;

\bigskip

\noindent $(\mathscr{A}_3):$ the sequence of decreasing eigenvalues  $\{\mu_q\}_{q\geq 1}$ of $\Sigma_Y$ is such that  $\sum_{q=1}^{+\infty}\mu_q^{1/2} <+\infty$;

\bigskip

\noindent $(\mathscr{A}_4):$ the operators $\Sigma_X$ and $\Sigma_Y$ have   infinitely many strictly positive eigenvalues; 

\bigskip

\noindent $(\mathscr{A}_5):$ the sequence $\{\gamma_n\}_{n\geq 1}$  satistfies: $\lim_{n\rightarrow +\infty}\Big(\gamma_n ^{-6} n^{-1}\Big)=0$. 

\bigskip

\noindent Then, we have:

\begin{theorem}\label{loi}
	Assume ($\mathscr{A}_1$)  to ($\mathscr{A}_5$).  Then,  under $\mathscr{H}_0$,  $\widehat{T}_n$ converges in distribution, as $n\rightarrow +\infty$,  to  $\mathcal{N}(0,1)$.
	
\end{theorem}
\bigskip

\noindent The resulting test for independence is performed as follows: for a given significance level $\alpha\in]0,1[$, one has to reject $\mathscr{H}_0$ if  $\vert \widehat{T}_n\vert> \mathbb{F}^{-1}(1-\alpha/2)$, where $\mathbb{F}$ is the cumulative distribution function of the standard normal distribution.
\subsection{Computational aspects}

\noindent For computing $\widehat{T}_n$  in practice, the kernel trick (see \cite{shawe04}) can be used. It is known from  \cite{fukumizu07b}  that $\widehat{N}_n  (X,Y)$  can be written as 
\begin{equation}\label{nchap}
\widehat{N}_n  (X,Y) 
 = \mbox{Tr} \bigg( G_X(G_X + n \gamma_n I_n )^{-1}G_Y(G_Y + n \gamma_n I_n )^{-1}\bigg),
\end{equation}
where $G_X$ and $G_Y$ are    the centered  $n\times n$ Gram matrices with terms equal to
$$ 
(G_X)_{ij} = \Big \langle K(X_i ,\cdot) - \overline{K}_n,K(X_j ,\cdot) - \overline{K}_n \Big \rangle_{\mathcal{H}_\mathcal{X}}\,\,\,
\textrm{ and }\,\,\,
(G_Y)_{ij} = \Big \langle L(Y_i ,\cdot) - \overline{L}_n,L(Y_j ,\cdot) - \overline{L}_n \Big \rangle_{\mathcal{H}_\mathcal{Y}},
$$
and $I_n$ is the $n\times n$ identity matrix. Clearly, using the reproducing properties of $K$ and $L$, the preceding terms can be writen as
\begin{equation*}\label{gram}
(G_X)_{ij} =k_{ij}-\frac{1}{n}\left(k_{i\centerdot}+k_{\centerdot j}\right)+\frac{k_{\centerdot\centerdot}}{n^2}\,\,\,
\textrm{ and }\,\,\,
(G_Y)_{ij} = \ell_{ij}-\frac{1}{n}\left(\ell_{i\centerdot}+\ell_{\centerdot j}\right)+\frac{\ell_{\centerdot\centerdot}}{n^2},
\end{equation*}
where
\[
k_{ij}=K(X_i,X_j),\,\,\,k_{i\centerdot}=\sum_{j=1}^nk_{ij},\,\,\,k_{\centerdot j}=\sum_{i=1}^nk_{ij},\,\,\,k_{\centerdot \centerdot}=\sum_{i=1}^n\sum_{j=1}^nk_{ij},
\]
and
\[
\ell_{ij}=L(Y_i,Y_j),\,\,\,\ell_{i\centerdot}=\sum_{j=1}^n\ell_{ij},\,\,\,\ell_{\centerdot j}=\sum_{i=1}^n\ell_{ij},\,\,\,\ell_{\centerdot \centerdot}=\sum_{i=1}^n\sum_{j=1}^n\ell_{ij}.
\]
Now let us look at the  computation of    $\widehat{D}_{1,n}$ and $\widehat{D}_{2,n}$. From the definition of $\widehat{\Sigma}_X$ it is easily seen that   any eigenvector    of this operator associated to a nonzero eigenvalue  is a linear combination of $u_1,\cdots,u_n$, where $ u_i = K(X_i ,\cdot) -\overline{K}_n$. So  $f= \sum_{j=1}^{n} \alpha_j u_j $   is an eigenvector of $\widehat \Sigma_{X}$ associated to the positive eigenvalue $\lambda$ if and only if
\[
  \frac{1}{n} \sum_{j=1}^{n} \alpha_j \sum_{i=1}^{n} \langle u_i , u_j \rangle_{\mathcal{H}_\mathcal{X}} u_i = \sum_{j=1}^{n} \lambda \alpha_j u_j,
\]
what is equivalent to
\[
\sum_{j=1}^{n} \Big \{ \frac{1}{n}\sum_{i=1}^{n} \alpha_i  (G_X)_{ij} \Big \} u_j = \sum_{j=1}^{n} \lambda \alpha_j u_j.
\]
This later equality is equivalent to 
\[
 \lambda \alpha_j = \frac{1}{n}\sum_{i=1}^{n} \alpha_i  (G_X)_{ij}, \; \; \; \forall \; 1 \leq j \leq n,
\]
that is  $ \frac{1}{n} G_X \alpha   =\lambda \alpha $.
 So the nonzero eigenvalues of $\widehat{\Sigma}_{X}$ are the same as these of   $\frac{1}{n} G_X$. Similarly,  the nonzero eigenvalues of  $\widehat{\Sigma}_{X} + \gamma_n \mathbb{I}_\mathcal{X}$ are the same as these of    $\frac{1}{n} (G_X + n \gamma_n I_n)$.  We then deduce that 
\[
d_{1}(\widehat{\Sigma}_X,\gamma_n)  = \mbox{Tr} (\widehat{\Sigma}_{X}(\widehat{\Sigma}_{X} + \gamma_n \mathbb{I}_\mathcal{X})^{-1})  = \mbox{Tr} \Big( G_X(G_X + n \gamma_n I_n )^{-1}\Big).
\]
From  similar reasoning we also get
\[
d_{1}(\widehat{\Sigma}_Y,\gamma_n)=\mbox{Tr} \Big( G_Y(G_Y + n \gamma_n I_n )^{-1}\Big),\,\,\,
d_{2}(\widehat{\Sigma}_X,\gamma_n)=\mbox{Tr} \Big( G_X^2(G_X + n \gamma_n I_n )^{-2}\Big)
\]
and
\[
d_{2}(\widehat{\Sigma}_Y,\gamma_n)=\mbox{Tr} \Big( G_Y^2(G_Y + n \gamma_n I_n )^{-2}\Big).
\]
Hence
\begin{equation}\label{calcd1}
\widehat{D}_{1,n} =   \mbox{Tr} \Big( G_X(G_X + n \gamma_n I_n )^{-1}\Big)\,\mbox{Tr} \Big( G_Y(G_Y + n \gamma_n I_n )^{-1}\Big)
\end{equation}
and
\begin{equation}\label{calcd2}
\widehat{D}_{2,n} = \bigg\{\mbox{Tr} \Big( G_X^2(G_X + n \gamma_n I_n )^{-2}\Big)\,\mbox{Tr} \Big( G_Y^2(G_Y + n \gamma_n I_n )^{-2}\Big)\bigg\}^{1/2}.
\end{equation}
Therefore $\widehat{T}_n$ is computed in practice by using \eqref{nchap}, \eqref{calcd1}  and \eqref{calcd2}. The crucial point is to compute $K(X_i,X_j)$ and $L(Y_i,Y_j)$ for all $(i,j)\in\{1,\cdots,n\}^2$.

\bigskip

\begin{remark}
In the case of  functional data corresponding, for instance,  to the case where the $X_i$s and the $Y_i$s are random functions belonging in $L^2([0,1])$ and observed on points $t_1,\cdots,t_r$ and $s_1,\cdots,s_q$, respectively, of   fine grids  in $[0,1]$ such that $t_1=s_1=0$ and $t_r=s_q=1$, the preceding terms can be computed or approximated easily, depending on   the used kernels. For example, 
if the gaussian kernel is used for $K$ and $L$ , one has
\[
K(X_i ,X_j )=\exp\left(-\omega^2\Vert X_i-X_j\Vert^2_\mathcal{X}\right)=\exp\left(-\omega^2\int_0^1\left(X_i(t)-X_j(t)\right)^2\,dt\right),
\]
and
\[
L(Y_i ,Y_j )=\exp\left(-\nu^2\Vert Y_i-Y_j\Vert^2_\mathcal{Y}\right)=\exp\left(-\nu^2\int_0^1\left(Y_i(t)-Y_j(t)\right)^2\,dt\right),
\]
where $\omega>0$ and $\nu>0$. These  terms  can    be approximated  by using trapezoidal rule, so leading to:
\begin{equation}\label{approx1}
K(X_i ,X_j )\simeq\exp\left(-\omega^2\sum_{m=1}^{r-1}\frac{t_{m+1}-t_m}{2}\left(\left(X_i(t_m)-X_j(t_m)\right)^2+\left(X_i(t_{m+1})-X_j(t_{m+1})\right)^2\right)\right)
\end{equation}
and
\begin{equation}\label{approx2}
L(Y_i ,Y_j )\simeq\exp\left(-\nu^2\sum_{m=1}^{q-1}\frac{s_{m+1}-s_m}{2}\left(\left(Y_i(s_m)-Y_j(s_m)\right)^2+\left(Y_i(s_{m+1})-Y_j(s_{m+1})\right)^2\right)\right).
\end{equation}
Then,   $\widehat{T}_{n}$ is  to be computed by using \eqref{approx1} and \eqref{approx2}.
\end{remark}

\section{Simulations}
\noindent In this section, we investigate the finite sample performances of the proposed test, that  we denote by PT,  and compare it to two  other kernel-based independence tests: the permutation test using HSIC  introduced in \cite{gretton05a}, denoted here by HSIC,  and the classical NOCCO-based permutation test of \cite{fukumizu07b}  that we denote by  pNOCCO. We computed   empirical sizes and powers    through   Monte Carlo simulations.  We considered  the case where $\mathcal{X}=\mathcal{Y}=\textrm{L}^2([0,1])$ and generated the data according to the two following models:

\bigskip

\noindent\textbf{Model 1:} $X(t) = t +\sqrt{2} \sum_{k=1}^{2}  \alpha_{k,1} \sin(k\pi t) + \alpha_{k,2}\cos(k\pi t)+ \varepsilon_1(t)  $ and $Y(t) = t + \sqrt{2}\sum_{k=1}^{2}  \beta_{k,1} \sin(k\pi t) + \beta_{k,2}  \cos(k\pi t)+ \varepsilon_2(t)  $, where $\varepsilon_1(t)$ and $\varepsilon_2(t)$ are independent and independently  sampled  from the normal idstribution $N(0,0.25)$, the Fourier  coefficients   are independently sampled from the standard normal distribution with the condition $\alpha_{2,2}=\beta_{2,2}$;

\bigskip

\noindent\textbf{Model 2:} $X(t)=\sqrt{2}\sum_{k=1}^{30}\xi_k\,\cos(k\pi t)$ and  $Y(t)=\sqrt{2}\sum_{k=1}^{30}\nu_k\,\cos(k\pi t)$, where the $\xi_k$s are independent and distributed as the Cauchy distribution $\mathscr{C}(0,0.5)$ and, for a given $m\in\{0,\cdots,30\}$, $\nu_k=f(\xi_k)$ for $k=1,\cdots,m$ and the $\nu_k$s with $k=m+1,\cdots,50$ are sampled independently from the standard normal distribution. Here, $f$ is a given function that establishes dependence between $X$ and $Y$.

\bigskip

\noindent In model 1, $X$ and $Y$ are dependent due to a shared coefficient. In model 2, $X$ and $Y$ are independent when $m=0$, and they are dependent for $m\in\{1,\cdots,30\}$, the dependence level increasing as $m$ increases. Empirical sizes and powers were computed on the basis of $500$ independent  replicates. For each of them, we generated a sample  of size   $n=20, 30, 50, 60, 100$ of the above  processes in discretized versions  on equispaced values $t_1,\cdots,t_{101}$  in  $[0,1]$, where $t_j=(j-1)/100$, $j=1,\cdots,101$.  For performing  our method, we used the gaussian  kernels $K(x,y)=L(x,y)=\exp\left(-\omega^2\int_0^1\left( x(t)-y(t)\right)^2dt\right)$ with bandwith $\omega^2$ equal to the heuristic median computed from the data;  it is the most popular bandwidth choice in simulations and while it has no guarantee of optimality, it remains a safe choice in most of the cases (see \cite{garreau18}).   The terms $K(X_i,X_j) $ and  $L(Y_i,Y_j) $
were computed by approximating   integrals involved in these   kernels  by using the trapezoidal rule as indicated in \eqref{approx1} and \eqref{approx2}.  The HSIC and pNOCCO methods    were  used with $100$ permutations. For PT and pNOCCO methods we used $\gamma_n=n^{-1/7}$ for the regularization sequence.  The nominal significance level is taken as $\alpha=0.05$ for all tests. Table 1 and Table 2  report  the obtained results.  Table 1 shows a superiority of  PT over  the two other methods. Indeed, the obtained values   for the power are higher for this method than for the two others for each sample size. In Table 2, the  obtained values for $m=0$  are closer to the nominal size for PT, and this method outperforms HSIC and pNOCCO for each values of $m$ greter than $0$. For all  the three methods the power increases as  $m$ increases, but it is seen that the power of PT increases faster to $1$.  In the particular case of $f(x)=cos(x)$, all the three  methods do not sufficiently detect the dependence since they give low power for all values of $m$.
\begin{figure}[! ht]
\centering
\begin{tabular}{ccccccc}
\hline
 $n$ & &\qquad PT & & HSIC & & pNOCCO \\
\hline
\hline
\\
 20&  &\qquad 0.504 & & 0.328 & & 0.354\\
\\
  30 & &\qquad 0.842 & & 0.586 & & 0.638\\
\\
  40 &  &\qquad 0.964 & & 0.784 & & 0.823\\
\\
  50 & &\qquad 0.992 & & 0.908 & & 0.942\\
\\
  60 & &\qquad 1.00 & & 0.966 & &  0.980\\
\\
\hline
\hline
\end{tabular}
\caption{Empirical powers over 500 replications for Model 1 with nominal significance level  $\alpha=0.05$.}
\end{figure}
\newpage
\begin{figure}[! ht]
\centering
\begin{tabular}{lclclclclcl}
\hline
Relationship & Method & $m=0$ & $m=1$  &  $m=3$ & $m=5$ & $m=10$ \\
\hline
\hline
&    &   &  & &  &   \\
 & PT & 0.030 & 0.184 & 0.558 & 0.858 & 1.000 \\
$f(x) = x^3$ &HSIC  & 0.010& 0.028 & 0.174 & 0.430 & 0.934 \\
 & pNOCCO & 0.008  & 0.036 & 0.192 & 0.450 & 0.942\\
&    &   &  & &  &   \\
\hline
&    &   &  & &  &   \\
 & PT & 0.023 & 0.170 & 0.530 & 0.844 & 0.996 \\
$f(x) = x^2$ & HSIC  & 0.010 & 0.032 & 0.164 & 0.432 & 0.926 \\
 & pNOCCO & 0.008 & 0.030 & 0.204 & 0.452 & 0.941\\
&    &   &  & &  &   \\
\hline
&    &   &  & &  &   \\
 & PT & 0.030 & 0.198 & 0.542 & 0.836 & 0.996 \\
$f(x) = x^2 sin(x)$ &HSIC  & 0.010 & 0.048 & 0.212 & 0.456 & 0.932 \\
 & pNOCCO & 0.008 & 0.048 & 0.218 & 0.462 & 0.936 \\
&    &   &  & &  &   \\
\hline
&    &   &  & &  &   \\
 &PT & 0.030 & 0.036 & 0.016 & 0.026 & 0.028 \\
$f(x) = sin(x)$ &HSIC  & 0.010 & 0.012 & 0.014 & 0.018 & 0.010 \\
 & pNOCCO & 0.008 & 0.014 & 0.012 & 0.016 & 0.018\\
&    &   &  & &  &   \\
\hline
\hline
\end{tabular}
\caption{Empirical sizes and powers  over 500 replications for Model 2, with sample size   $n=100$ and   significance level  $\alpha=0.05$.}
\end{figure}

\section{Proofs}

\subsection{Preliminary   lemmas}
\noindent  First, putting
\[
\widetilde{\Sigma}_{XY}   = \frac{1}{n} \sum_{i=1}^{n}  (L(Y_j,) - m_Y) \otimes (K(X_j,.) - m_X),\,\,\, \widetilde{V}_{XY}    = (\Sigma_{X} + \gamma_n I)^{-1/2}\widehat{\Sigma}_{XY}   ( \Sigma_{Y} + \gamma_n I)^{-1/2}\,\,\,
\]
and
\[
\widetilde{N}_n  (X,Y) = \Vert\widetilde{V}_{XY}  \Vert_{\textrm{HS}} ^2,
\]
we have:

\bigskip

\begin{lemma}\label{lemme1}
Assume ($\mathscr{A}_1$) to ($\mathscr{A}_5$). Then, under $\mathscr{H}_0$, we have:
\[
\big\vert \widehat{N}_n (X,Y) -  \widetilde{N}_n  (X,Y) \big\vert = O_P (\gamma_n ^{-3} n^{-3/2}),\,\,\,
  \widehat{N}_n (X,Y) = O_P (\gamma_n ^{-2} n^{-1})
\]   
and    
\[
\widetilde{N}_n  (X,Y)  = O_P (\gamma_n ^{-2} n^{-1}).
\]
\end{lemma}
\noindent\textit{Proof}.
We have
\begin{equation}\label{Nchap}
\begin{aligned}
\widehat{N}_n (X,Y) &   =  \Big \langle \left(\widehat{\Sigma}_{X} + \gamma_n \mathbb{I}_\mathcal{X}\right)^{-1}\widehat{\Sigma}_{XY}   \; , \; \widehat{\Sigma}_{XY}  \left( \widehat{\Sigma}_{Y} + \gamma_n \mathbb{I}_\mathcal{Y}\right)^{-1} \Big  \rangle_{\textrm{HS}}\\
\end{aligned}
\end{equation}
and, similarly,
\begin{equation}\label{Ntilde}
 \widetilde{N}_n  (X,Y)=\Big \langle \left( \Sigma_{X} + \gamma_n \mathbb{I}_\mathcal{X}\right)^{-1}\widehat{\Sigma}_{XY}   \; , \; \widehat{\Sigma}_{XY}  \left(  \Sigma_{Y} + \gamma_n \mathbb{I}_\mathcal{Y}\right)^{-1} \Big  \rangle_{\textrm{HS}},
\end{equation}
so that
\begin{eqnarray*}
& &\widehat{N}_n (X,Y) -  \widetilde{N}_n  (X,Y) \\
& = &\Big \langle\Big(\left(\widehat{\Sigma}_{X} + \gamma_n \mathbb{I}_\mathcal{X}\right)^{-1}   - \left( \Sigma_{X} + \gamma_n \mathbb{I}_\mathcal{X}\right)^{-1}\Big)\widehat{\Sigma}_{XY}   \; , \; \widehat{\Sigma}_{XY}    \left( \widehat{\Sigma}_{Y} + \gamma_n \mathbb{I}_\mathcal{Y}\right)^{-1} \Big  \rangle_{\textrm{HS}}\\
&  &\qquad \qquad + \Big \langle\left( \Sigma_{X} + \gamma_n \mathbb{I}_\mathcal{X}\right)^{-1}\widehat{\Sigma}_{XY}   \; , \; \widehat{\Sigma}_{XY} \Big( \left( \widehat{\Sigma}_{Y} + \gamma_n \mathbb{I}_\mathcal{Y}\right)^{-1} -  \left(  \Sigma_{Y} + \gamma_n \mathbb{I}_\mathcal{Y}\right)^{-1}\Big) \Big  \rangle_{\textrm{HS}} \\
&  &:=A_1 +A_2.
\end{eqnarray*}
Using the formula $A^{-1} - B^{-1} = A^{-1}(B-A)B^{-1}$ we get
\begin{equation*}
\begin{aligned}
\Big(\left(\widehat{\Sigma}_{X} + \gamma_n \mathbb{I}_\mathcal{X}\right)^{-1}   - \left( \Sigma_{X} + \gamma_n \mathbb{I}_\mathcal{X}\right)^{-1}\Big)\widehat{\Sigma}_{XY}  
& =\left(\widehat{\Sigma}_{X} + \gamma_n \mathbb{I}_\mathcal{X}\right)^{-1}\Big (  \Sigma_{X}  -  \widehat{\Sigma}_{X} \Big) \left( \Sigma_{X} + \gamma_n \mathbb{I}_\mathcal{X}\right)^{-1}\Big)\widehat{\Sigma}_{XY}.  \\
\end{aligned}
\end{equation*}
Then, from the  Cauchy-Schwartz inequality and  $\Vert( \Sigma + \gamma_n I)^{-1}\Vert \leq \gamma_n ^{-1}$ for all compact operator $\Sigma$, we get  
\begin{equation}\label{i}
\begin{aligned}
|A_1| & \leq \gamma_n ^{-3}  \Vert\Sigma_{X}  -  \widehat{\Sigma}_{X}\Vert_{\textrm{HS}}\Vert\widehat{\Sigma}_{XY}  \Vert^2_{\textrm{HS}} 
\end{aligned}
\end{equation}
and, similarly, 
\begin{equation}\label{ii}
\begin{aligned}
|A_2| & \leq \gamma_n ^{-3}   \Vert\Sigma_{Y}  -  \widehat{\Sigma}_{Y}\Vert_{\textrm{HS}}\Vert\widehat{\Sigma}_{XY}  \Vert_{\textrm{HS}}^2. 
\end{aligned}
\end{equation}
Lemma 5 in \cite{fukumizu07a}  ensures that $\Vert\widehat{\Sigma}_{XY} - {\Sigma}_{XY}\Vert_{\textrm{HS}}=O_P (n^{-1/2})$, $\Vert\Sigma_{Y}  -  \widehat{\Sigma}_{Y}\Vert_{\textrm{HS}}=O_P (n^{-1/2})$ and $\Vert\Sigma_{X}  -  \widehat{\Sigma}_{X}\Vert_{\textrm{HS}}=O_P (n^{-1/2})$. Under $\mathscr{H}_0$ one has $\Sigma_{XY}=0$, hence   $\Vert\widehat{\Sigma}_{XY}\Vert_{\textrm{HS}}=O_P (n^{-1/2})$. We the deduce from \eqref{i} and \eqref{ii} that   $\vert A_1\vert =O_P(\gamma_n ^{-3} n^{-3/2})$, and  $\vert A_2\vert =O_P(\gamma_n ^{-3} n^{-3/2})$. Consequently,  $\big| \widehat{N}_n (X,Y) -  \widetilde{N}_n  (X,Y) \big|  =O_P(\gamma_n ^{-3} n^{-3/2})$. Now, from   \eqref{Nchap} and \eqref{Ntilde}  we get by the Cauchy-Schwartz inequality
$ |\widehat{N}_n (X,Y)| \leq \gamma_n ^{-2} \Vert\widehat{\Sigma}_{XY}\Vert_{\textrm{HS}} ^2$
and
$|\widetilde{N}_n (X,Y)| \leq \gamma_n ^{-2} \Vert\widehat{\Sigma}_{XY}\Vert_{\textrm{HS}} ^2$
ensuring that $\widehat{N}_n (X,Y) = O_P (\gamma_n ^{-2} n^{-1})$ and $\widetilde{N}_n (X,Y) = O_P (\gamma_n ^{-2} n^{-1})$.

\bigskip

\noindent Secondly, for $r\in\{1,2\}$, we consider  
\[
 D_{r,n}=  d_{r}(\Sigma_X,\gamma_n)d_{r}(\Sigma_Y,\gamma_n)=\bigg(     \sum_{p=1}^{+\infty}   \frac{\lambda_p^r}{\left(\lambda_p + \gamma_n\right)^r}     \bigg)^{1/r}   \bigg(     \sum_{p=1}^{+\infty}   \frac{\mu_p^r}{\left(\mu_p + \mu_n\right)^r}     \bigg)^{1/r},
\]
where  $\{\lambda_p\}_{p\geq 1}$ (resp. $\{\mu_q\}_{q\geq 1}$ )   is the sequence of decreasing eigenvalues  of $\Sigma_X$ (resp. $\Sigma_Y$), and we put:
\[
\widetilde{T}_n =  \frac{n  \widetilde{N}_n  (X,Y) -   D_{1,n}}{\sqrt{2} D_{2,n}}.  
\]
Then, we have:

\bigskip

\begin{lemma}\label{lem2}
Assume ($\mathscr{A}_1$) to ($\mathscr{A}_5$). Then, under $\mathscr{H}_0$, we have $\widehat{T_n}=\widetilde{T}_n+ o_p(1)$.
\end{lemma}
\noindent\textit{Proof}.
We have
\begin{equation}\label{Tn}
\begin{aligned}
\sqrt{2}\big( \widehat{T}_n  - \widetilde{T}_n\big) = \bigg( \frac{n \widehat{N}_n (X,Y) }{\widehat{D}_{2,n}} - \frac{n  \widetilde{N}_n  (X,Y)}{ D_{2,n}}\bigg) - \bigg( \frac{\widehat{D}_{1,n}}{\widehat{D}_{2,n}} - \frac{ D_{1,n}}{ D_{2,n}}\bigg): = B- C.
\end{aligned}
\end{equation}
First,
\begin{equation}\label{A}
\begin{aligned}
\vert B\vert & \leq \bigg\vert \frac{1}{\widehat{D}_{2,n}} - \frac{1}{ D_{2,n}}\bigg|n \widehat{N}_n (X,Y) + \frac{n}{ D_{2,n}} \bigg\vert \widehat{N}_n (X,Y) -  \widetilde{N}_n  (X,Y) \bigg\vert\\
& \leq  \bigg\vert \frac{ D_{2,n} - \widehat{D}_{2,n}}{ D_{2,n}\widehat{D}_{2,n}}\bigg\vert \; \big\vert n\widehat{N}_n (X,Y) - n \widetilde{N}_n  (X,Y) \big| \\
& +  \bigg\vert \frac{ D_{2,n} - \widehat{D}_{2,n}}{ D_{2,n}\widehat{D}_{2,n}}\bigg\vert n \widetilde{N}_n  (X,Y) +   \frac{1}{ D_{2,n}} \bigg\vert n\widehat{N}_n (X,Y) -  n\widetilde{N}_n  (X,Y) \bigg\vert \\
& := B_1 + B_2 + B_3.
\end{aligned}
\end{equation}
Clearly,
\begin{equation*}
\begin{aligned}
\frac{ D_{2,n} - \widehat{D}_{2,n}}{ D_{2,n}\widehat{D}_{2,n}} 
& = \frac{ d_{2}({\Sigma}_X,\gamma_n) - d_{2}(\widehat{\Sigma}_X,\gamma_n) }{d_{2}({\Sigma}_X,\gamma_n)d_{2}(\widehat{\Sigma}_X,\gamma_n)d_{2}(\widehat{\Sigma}_Y,\gamma_n)} + \frac{ d_{2}({\Sigma}_Y,\gamma_n)-d_{2}(\widehat{\Sigma}_Y,\gamma_n) }{d_{2}({\Sigma}_X,\gamma_n)d_{2}({\Sigma}_Y,\gamma_n)d_{2}(\widehat{\Sigma}_Y,\gamma_n)}
\end{aligned}
\end{equation*}
and, therefore,
\begin{equation*}\label{Achap2}
\begin{aligned}
\vert B_1\vert  & \leq  \frac{\big \vert d_{2}({\Sigma}_X,\gamma_n) - d_{2}(\widehat{\Sigma}_X,\gamma_n)  \big \vert}{d_{2}({\Sigma}_X,\gamma_n)d_{2}(\widehat{\Sigma}_X,\gamma_n)d_{2}(\widehat{\Sigma}_Y,\gamma_n)}\times \big| n\widehat{N}_n (X,Y) - n \widetilde{N}_n  (X,Y) \big|  \\
& \qquad \qquad  \qquad + \frac{\big\vert d_{2}({\Sigma}_Y,\gamma_n)-d_{2}(\widehat{\Sigma}_Y,\gamma_n)\big\vert}{d_{2}({\Sigma}_X,\gamma_n)d_{2}({\Sigma}_Y,\gamma_n)d_{2}(\widehat{\Sigma}_Y,\gamma_n)} \times \big| n\widehat{N}_n (X,Y) - n \widetilde{N}_n  (X,Y) \big\vert ,\\
\vert B_2\vert & \leq  \frac{\big \vert d_{2}({\Sigma}_X,\gamma_n) - d_{2}(\widehat{\Sigma}_X,\gamma_n)\big\vert}{d_{2}({\Sigma}_X,\gamma_n)d_{2}(\widehat{\Sigma}_X,\gamma_n)d_{2}(\widehat{\Sigma}_Y,\gamma_n)}\times n \widetilde{N}_n  (X,Y) \\
& \qquad \qquad  \qquad + \frac{\vert d_{2}({\Sigma}_Y,\gamma_n)-d_{2}(\widehat{\Sigma}_Y,\gamma_n)\big\vert}{d_{2}({\Sigma}_X,\gamma_n)d_{2}({\Sigma}_Y,\gamma_n)d_{2}(\widehat{\Sigma}_Y,\gamma_n)} \times n \widetilde{N}_n  (X,Y).
\end{aligned}
\end{equation*}
Proposition 12  in \cite{harchaoui08}  and Assumption   ($\mathscr{A}_5$)  ensure that 
$$\Vert\widehat{\Sigma}_X - \Sigma_X\Vert_{\mathcal C_1}:=\sum_{p=1}^{+\infty} \Vert(\widehat{\Sigma}_X - \Sigma_X)e_p\Vert_{\mathcal{H}_{\mathcal{X}}} = O_P (n^{-1/2});$$ 
applying now Lemma 23 in  \cite{harchaoui08}   with $S=\Sigma_X$ and $\Delta = \widehat{\Sigma}_X - \Sigma_X$ we get
\begin{equation}\label{d2}
\begin{aligned}
\big\vert d_{2}({\Sigma}_X,\gamma_n) - d_{2}(\widehat{\Sigma}_X,\gamma_n)  \big\vert   \leq \frac{\gamma_n ^{-1} \Vert \widehat{\Sigma}_X - \Sigma_X \Vert_{\mathcal C_1}}{1 + \gamma_n ^{-1} \Vert \widehat{\Sigma}_X - \Sigma_X \Vert_{\mathcal C_1}} = O_P (\gamma_n ^{-1} n^{-1/2})
\end{aligned}
\end{equation}
and, similarly, $
\big\vert d_{2}({\Sigma}_Y,\gamma_n) - d_{2}(\widehat{\Sigma}_Y,\gamma_n)  \big\vert  = O_P (\gamma_n ^{-1} n^{-1/2})$, 
$\big\vert d_{1}({\Sigma}_X,\gamma_n) - d_{1}(\widehat{\Sigma}_X,\gamma_n)  \big\vert  = O_P (\gamma_n ^{-1} n^{-1/2})$ and 
$\big\vert d_{1}({\Sigma}_Y,\gamma_n) - d_{1}(\widehat{\Sigma}_Y,\gamma_n)  \big\vert  = O_P (\gamma_n ^{-1} n^{-1/2})$. Plugging these equalities in    (\ref{d2}) and using Lemma 1 and Lemma 18 in  \cite{harchaoui08}, we deduce from \eqref{A}  that $B = o_P(1)$. On the other hand, 
\begin{equation*}
\begin{aligned}
C & = \frac{d_{1}(\widehat{\Sigma}_X,\gamma_n)d_{1}(\widehat{\Sigma}_Y,\gamma_n)}{d_{2}(\widehat{\Sigma}_X,\gamma_n)d_{2}(\widehat{\Sigma}_Y,\gamma_n)} - \frac{d_{1}({\Sigma}_X,\gamma_n)d_{1}({\Sigma}_Y,\gamma_n)}{d_{2}({\Sigma}_X,\gamma_n)d_{2}({\Sigma}_Y,\gamma_n)}\\
& = \frac{\big( d_{1}(\widehat{\Sigma}_X,\gamma_n) - d_{1}(\widehat{\Sigma}_X,\gamma_n) \big )d_{1}(\widehat{\Sigma}_Y,\gamma_n)}{d_{2}(\widehat{\Sigma}_X,\gamma_n)d_{2}(\widehat{\Sigma}_Y,\gamma_n)} + \frac{d_{1}(\widehat{\Sigma}_X,\gamma_n)\big ( d_{1}(\widehat{\Sigma}_Y,\gamma_n)  -  d_{1}({\Sigma}_Y,\gamma_n) \big)}{d_{2}(\widehat{\Sigma}_X,\gamma_n)d_{2}(\widehat{\Sigma}_Y,\gamma_n)} \\
&   + \frac{d_{1}({\Sigma}_X,\gamma_n)d_{1}({\Sigma}_Y,\gamma_n) \big( d_{2}(\widehat{\Sigma}_X,\gamma_n) - d_{2}({\Sigma}_X,\gamma_n) \big )}{d_{2}(\widehat{\Sigma}_X,\gamma_n)d_{2}({\Sigma}_X,\gamma_n)d_{2}({\Sigma}_Y,\gamma_n)} \\
&+ \frac{d_{1}({\Sigma}_X,\gamma_n)d_{1}({\Sigma}_Y,\gamma_n)\big ( d_{2}(\widehat{\Sigma}_Y,\gamma_n) - d_{2}({\Sigma}_Y,\gamma_n) \big )}{d_{2}(\widehat{\Sigma}_Y,\gamma_n)d_{2}({\Sigma}_X,\gamma_n)d_{2}({\Sigma}_Y,\gamma_n)};
\end{aligned}
\end{equation*}
so that 
\begin{equation*}
\begin{aligned}
\vert C\vert  &\leq \frac{O_P (\gamma_n ^{-1} n^{-1/2})d_{1}(\widehat{\Sigma}_Y,\gamma_n)}{d_{2}(\widehat{\Sigma}_X,\gamma_n)d_{2}(\widehat{\Sigma}_Y,\gamma_n)} + \frac{O_P (\gamma_n ^{-1} n^{-1/2})d_{1}(\widehat{\Sigma}_X,\gamma_n)}{d_{2}(\widehat{\Sigma}_X,\gamma_n)d_{2}(\widehat{\Sigma}_Y,\gamma_n)} \\
&+ \frac{O_P (\gamma_n ^{-1} n^{-1/2})d_{1}({\Sigma}_X,\gamma_n)d_{1}({\Sigma}_Y,\gamma_n)}{d_{2}(\widehat{\Sigma}_X,\gamma_n)d_{2}({\Sigma}_X,\gamma_n)d_{2}({\Sigma}_Y,\gamma_n)}\\
& + \frac{O_P (\gamma_n ^{-1} n^{-1/2})d_{1}({\Sigma}_X,\gamma_n)d_{1}({\Sigma}_Y,\gamma_n)}{d_{2}(\widehat{\Sigma}_Y,\gamma_n)d_{2}({\Sigma}_X,\gamma_n)d_{2}({\Sigma}_Y,\gamma_n)}\\
& = \frac{O_P (\gamma_n ^{-2} n^{-1/2}) \gamma_n d_{1}(\widehat{\Sigma}_Y,\gamma_n)}{d_{2}(\widehat{\Sigma}_X,\gamma_n)d_{2}(\widehat{\Sigma}_Y,\gamma_n)}  + \frac{O_P (\gamma_n ^{-2} n^{-1/2}) \gamma_n d_{1}(\widehat{\Sigma}_X,\gamma_n)}{d_{2}(\widehat{\Sigma}_X,\gamma_n)d_{2}(\widehat{\Sigma}_Y,\gamma_n)} \\
&+ \frac{O_P (\gamma_n ^{-3} n^{-1/2}) \gamma_n d_{1}({\Sigma}_X,\gamma_n) \gamma_n d_{1}({\Sigma}_Y,\gamma_n)}{d_{2}(\widehat{\Sigma}_X,\gamma_n)d_{2}({\Sigma}_X,\gamma_n)d_{2}({\Sigma}_Y,\gamma_n)}\\
& + \frac{O_P (\gamma_n ^{-3} n^{-1/2}) \gamma_n d_{1}({\Sigma}_X,\gamma_n) \gamma_n d_{1}({\Sigma}_Y,\gamma_n)}{d_{2}(\widehat{\Sigma}_Y,\gamma_n)d_{2}({\Sigma}_X,\gamma_n)d_{2}({\Sigma}_Y,\gamma_n)}\\
\end{aligned}
\end{equation*}
Lemmas 18 and 19 in \cite{harchaoui08}   and  Assumption   ($\mathscr{A}_5$) allow to conclude that $C=o_P(1)$. Then from \eqref{Tn} we deduce that, under $\mathscr{H}_0$, $\widehat{T}_n  = \widetilde{T}_n+ o_p(1)$.

\bigskip

\noindent Next,  we can write $\widehat{\Sigma}_{XY} $ as
$$ \widehat{\Sigma}_{XY}   = \frac{1}{n} \sum_{i=1}^{n}  \Big(L(Y_i,\cdot)-m_Y\Big) \otimes \Big(K(X_i,\cdot)-m_X\Big) - \Big (\overline{L}_n - m_Y\Big) \otimes \Big(\overline{K}_n -m_X\Big),  $$
where
\[
\overline{L}_n=\frac{1}{n} \sum_{i=1}^{n} L(Y_i,\cdot)\,\,\,\textrm{ and }\,\,\,\overline{K}_n=\frac{1}{n} \sum_{i=1}^{n} K(X_i,\cdot).
\]
Let $\{e_p\}_{p\geq 1}$ (resp. $\{f_q\}_{q\geq 1}$) be an orthonomal basis of $\mathcal{H}_\mathcal{X}$ (resp. $\mathcal{H}_\mathcal{Y}$) consisting of eigenvectors of $\Sigma_X$ (resp. $\Sigma_Y$) and being such that $e_p$ (resp. $f_q$) is eigenvector associated with the eigenvalue $\lambda_p$ (resp. $\mu_q$). We have the the eigen-decompositions
\begin{equation*}
(\Sigma_{X} + \gamma_n \mathbb{I}_\mathcal{X})^{-1/2} = \sum_{p=1}^{+\infty} {(\lambda_p + \gamma_n)}^{-1/2} {e}_p\otimes {e}_p,\,\,\,(\Sigma_{Y} + \mathbb{I}_\mathcal{Y})^{-1/2} = \sum_{q=1}^{+\infty} {(\mu_q + \gamma_n)}^{-1/2} f_q\otimes f_q
\end{equation*}
lead to $ \widetilde{V}_{XY}=\widetilde{V}_{XY} ^{(1)} - \widetilde{V}_{XY} ^{(2)}$,
where
\[
\widetilde{V}_{XY} ^{(1)}=\frac{1}{n} \sum_{p=1}^{+\infty} \sum_{q=1}^{+\infty}\sum_{i=1}^{n} (\lambda_p + \gamma_n)^{-1/2}(\mu_q + \gamma_n)^{-1/2} (e_p \otimes e_p) \bigg( (L(Y_i,\cdot)-m_Y)\otimes ( K(X_i ,\cdot)- m_X) \bigg)(f_q \otimes f_q)
\]
and
\[
\widetilde{V}_{XY} ^{(2)}=\sum_{p=1}^{+\infty} \sum_{q=1}^{+\infty} {(\lambda_p + \gamma_n)}^{-1/2}{(\mu_q + \gamma_n)}^{-1/2} (e_p \otimes e_p)\big( ( \overline{L}_n - m_Y )  \otimes (\overline{K}_n -m_X)  \big)(f_q \otimes f_q).
\]
The following lemma gives properties of the last two terms that will be useful later for proving the main theorem of the paper.

\bigskip

\begin{lemma}\label{lem3}
Assume ($\mathscr{A}_1$) to ($\mathscr{A}_5$). Then, under $\mathscr{H}_0$, we have:
\begin{itemize}
\item[(i)] $\mathbb{E} \Big (n\Vert \widetilde{V}_{XY}^{(1)} \Vert_{\textrm{HS}} ^2\Big )= D_{1,n}$;
\item[(ii)]  $\mathbb{E} \Big (n\Vert \widetilde{V}_{XY}^{(2)} \Vert_{\textrm{HS}} ^2\Big )=n^{-1} D_{1,n}$;
\item[(iii)]
 $n  \Vert  \widetilde{V}_{XY} ^{(2)}\Vert_{\textrm{HS}}^2  =o_P (1)$;
\item[(iv)]$n \big \langle  \widetilde{V}_{XY} ^{(1)} \; , \widetilde{V}_{XY} ^{(2)} \big \rangle_{\textrm{HS}} = o_P (1)$.
\end{itemize}
\end{lemma}
\noindent\textit{Proof}.

\noindent\textit{(i)}. Using twice the formula $(a \otimes b)(c \otimes d) = \langle a , d \rangle c \otimes b$ and the reproducing properties of $K$ and $L$, we get
\begin{equation*}
\begin{aligned}
\widetilde{V}_{XY}^{(1)}   
& =\frac{1}{n} \sum_{p=1}^{+\infty} \sum_{q=1}^{+\infty} {(\lambda_p + \gamma_n)}^{-1/2}{(\mu_q + \gamma_n)}^{-1/2}\\
&\hspace{2cm} \bigg(\sum_{i=1}^{n}\big( e_p (X_i) - \mathbb{E} (e_p (X_i))  \big( f_q (Y_i) - \mathbb{E}(f_q (Y_i)) \Big) \bigg)f_q \otimes e_p;\\
\end{aligned}
\end{equation*}
so that
\begin{equation}\label{helpfull1}
\begin{aligned}
\Vert \widetilde{V}_{XY} ^{(1)}\Vert_{\textrm{HS}} ^2 &=  \frac{1}{n^2} \sum_{p=1}^{+\infty} \sum_{q=1}^{+\infty}  {(\lambda_p + \gamma_n)}^{-1}{(\mu_q + \gamma_n)}^{-1}   \Big  \{  \sum_{i=1}^{n} \alpha_{p,i,q}\Big \}^2\\
& =  \frac{1}{n^2} \sum_{p=1}^{+\infty} \sum_{q=1}^{+\infty}  {(\lambda_p + \gamma_n)}^{-1}{(\mu_q + \gamma_n)}^{-1} \bigg(  \sum_{i=1}^{n} \alpha_{p,i,q}^2 + 2  \sum_{i=2}^{n} \alpha_{p,i,q}M_{p,i-1,q} \bigg),\\
\end{aligned}
\end{equation}
where 
\begin{equation}\label{alpha}
\alpha_{p,i,q} = \big( e_p (X_i) - \mathbb{E} (e_p (X_i)) \big) \big( f_q (Y_i) - \mathbb{E}(f_q (Y_i)) \big)
\end{equation}
and
\begin{equation}\label{mpjq}
M_{p,j,q}= \sum_{l=1}^{j} \alpha_{p,l,q}.
\end{equation}
Since  $\alpha_{p,i,q}$ and $\alpha_{p,j,q}$   are independent  for $i\neq j$  and,  under $\mathscr{H}_0$, $\mathbb{E}\left(\alpha_{p,i,q}\right)=0$,  we get
\begin{equation}\label{helpfull2}
\begin{aligned}
\mathbb{E} \bigg(n  \Vert \widetilde{V}_{XY} ^{(1)}\Vert_{\textrm{HS}} ^2\bigg)& =  \frac{1}{n} \sum_{p=1}^{+\infty} \sum_{q=1}^{+\infty}  {(\lambda_p + \gamma_n)}^{-1}{(\mu_q + \gamma_n)}^{-1}\\
&\hspace{2cm}  \Bigg(  \sum_{i=1}^{n} \mathbb{E} \bigg(\alpha_{p,i,q}^2\bigg) + 2  \sum_{i=2}^{n}  \mathbb{E} \bigg(\alpha_{p,i,q}M_{p,i-1;q}\bigg) \Bigg)  \\
& =   \frac{1}{n} \sum_{p=1}^{+\infty} \sum_{q=1}^{+\infty}  {(\lambda_p + \gamma_n)}^{-1}{(\mu_q + \gamma_n)}^{-1} \\
&\hspace{2cm}   \Bigg(  \sum_{i=1}^{n} \mathbb{E} \bigg(\alpha_{p,i,q}^2\bigg) + 2  \sum_{i=2}^{n}  \mathbb{E} \bigg(\alpha_{p,i,q}\bigg)\,\mathbb{E}\bigg(M_{p,i-1;q}\bigg) \Bigg) \\
& =   \frac{1}{n} \sum_{p=1}^{+\infty} \sum_{q=1}^{+\infty}  {(\lambda_p + \gamma_n)}^{-1}{(\mu_q + \gamma_n)}^{-1}  \Bigg(  \sum_{i=1}^{n} \mathbb{E} \bigg(\alpha_{p,i,q}^2\bigg)  \Bigg) \\
& =  \sum_{p=1}^{+\infty} \sum_{q=1}^{+\infty}  {(\lambda_p + \gamma_n)}^{-1}{(\mu_q + \gamma_n)}^{-1}  var (e_p (X_1)) var (f_q(Y_1)) \\
& =  \sum_{p=1}^{+\infty} \sum_{q=1}^{+\infty}  {(\lambda_p + \gamma_n)}^{-1}{(\mu_q + \gamma_n)}^{-1}  \langle e_p , \Sigma_X e_p \rangle_{\mathcal{H}_{\mathcal{X}}} \langle f_q , \Sigma_Y f_q \rangle_{\mathcal{H}_{\mathcal{Y}}} \\
& =  \sum_{p=1}^{+\infty} \sum_{q=1}^{+\infty}  {(\lambda_p + \gamma_n)}^{-1}{(\mu_q + \gamma_n)}^{-1} \lambda_p \mu_q    \\
& = \Big \{     \sum_{p=1}^{+\infty}   \frac{\lambda_p}{\lambda_p + \gamma_n}     \Big \}   \Big \{ \sum_{q=1}^{+\infty}     \frac{\mu_q}{\mu_q + \gamma_n}  \Big \}     \\ 
 & = d_{1}(\Sigma_X,\gamma_n)d_{1}(\Sigma_Y,\gamma_n)= D_{1,n}.
\end{aligned}
\end{equation}
\noindent\textit{(ii)}.  Using again  the formula $(a \otimes b)(c \otimes d) = \langle a , d \rangle\, c \otimes b$ and the reproducing properties of $K$ and $L$, we get
\begin{equation*}
\begin{aligned}
\widetilde{V}_{XY}^{(2)}   & =\frac{1}{n^2} \sum_{p=1}^{+\infty} \sum_{q=1}^{+\infty} {(\lambda_p + \gamma_n)}^{-1/2}{(\mu_q + \gamma_n)}^{-1/2}\\
&\hspace{2cm}   \Big( \sum_{i=1}^{n}e_p (X_i) - \mathbb{E}(e_p (X_i)) \Big) \Big( \sum_{i=1}^{n}f_q (Y_i) - \mathbb{E} (f_q (Y_i)) \Big)f_q \otimes e_p,\\
\end{aligned}
\end{equation*}
so that 
\begin{equation*}
\begin{aligned}
\Vert \widetilde{V}_{XY}^{(2)} \Vert_{\textrm{HS}} ^2  & =\frac{1}{n^4} \sum_{p=1}^{+\infty} \sum_{q=1}^{+\infty} {(\lambda_p + \gamma_n)}^{-1}{(\mu_q + \gamma_n)}^{-1}\\
&\hspace{2cm}   \Big( \sum_{i=1}^{n}e_p (X_i) - \mathbb{E} (e_p (X_i)) \Big)^2 \Big( \sum_{i=1}^{n}f_q (Y_i) - \mathbb{E} (f_q (Y_i)) \Big)^2\\
\end{aligned}
\end{equation*}
and, under $\mathscr{H}_0$,
\begin{equation*}
\begin{aligned}
\mathbb{E} \Big (n\Vert \widetilde{V}_{XY}^{(2)} \Vert_{\textrm{HS}} ^2\Big )  & =\frac{1}{n^3} \sum_{p=1}^{+\infty} \sum_{q=1}^{+\infty} {(\lambda_p + \gamma_n)}^{-1}{(\mu_q + \gamma_n)}^{-1} \\
&\hspace{2cm}  \mathbb{E} \Bigg (\bigg( \sum_{i=1}^{n}e_p (X_i) - \mathbb{E} (e_p (X_i)) \bigg)^2 \Bigg)\, \mathbb{E} \Bigg (\bigg( \sum_{i=1}^{n}f_q (Y_i) - \mathbb{E} (f_q (Y_i)) \bigg)^2 \Bigg )\\
& = \frac{1}{n^3} \sum_{p=1}^{+\infty} \sum_{q=1}^{+\infty} {(\lambda_p + \gamma_n)}^{-1}{(\mu_q + \gamma_n)}^{-1}  \Bigg(\sum_{i=1}^{n}\mbox{var}\,\Big(e_p (X_i)\Big)\Bigg)\,\Bigg(\sum_{j=1}^{n} \mbox{var}\,\Big( f_q (Y_j)\Big)\Bigg)\\
& = \frac{1}{n}\sum_{p=1}^{+\infty} \sum_{q=1}^{+\infty} {(\lambda_p + \gamma_n)}^{-1}{(\mu_q + \gamma_n)}^{-1}  \lambda_p \mu_q\\
& = \frac{1}{n} D_{1,n}.
\end{aligned}
\end{equation*}
\noindent\textit{(iii)}. Lemma 19  in \cite{harchaoui08}   and Assumptions ($\mathscr{A}_2$) and  ($\mathscr{A}_3$)  imply that 
\begin{equation*}
\begin{aligned}
\gamma_n  D_{1,n} & \leq 4  \bigg(\sum_{ p \geq 1} \lambda_p ^{1/2}\bigg)\bigg(\sum_{ q \geq 1} \mu_q ^{1/2}\bigg).
\end{aligned}
\end{equation*}
Then,  using Assumption ($\mathscr{A}_5$), we obtain the inequality
\begin{equation}\label{ineg}
\begin{aligned}
\mathbb{E} \Big (n\Vert \widetilde{V}_{XY}^{(2)} \Vert_{\textrm{HS}} ^2\Big )=n^{-1}  D_{1,n}  
& \leq 4 \gamma_n ^{-1}n^{-1} \bigg(\sum_{ p \geq 1} \lambda_p ^{1/2}\bigg)\bigg(\sum_{ q \geq 1} \mu_q ^{1/2}\bigg) = o(1)
\end{aligned}
\end{equation}
from which, and Markov inequality, we deduce that    $n\Vert \widetilde{V}_{XY}^{(2)} \Vert_{\textrm{HS}} ^2=o_P(1)$ since $\lim_{n\rightarrow +\infty}\left(\gamma_n ^{-1}n^{-1}\right)=0$.

\noindent\textit{(iv)}. Using the  Cauchy-Schwartz and Hölder inequalities, the previous properties $(i)$ and $(ii)$, and \eqref{ineg}, we obtain the inequality 
\begin{equation*}
\begin{aligned}
\mathbb{E} \Bigg(n\big\vert \big  \langle  \widetilde{V}_{XY} ^{(1)} \; , \widetilde{V}_{XY} ^{(2)} \big \rangle_{\textrm{HS}}\bigg\vert  \Bigg) & \leq n\Bigg(\mathbb{E} \big\Vert \widetilde{V}_{XY} ^{(1)} \bigg\vert  \bigg\vert  _{\textrm{HS}} ^2\Bigg)^{1/2}\Bigg(\mathbb{E} \big\Vert \widetilde{V}_{XY} ^{(2)} \bigg\vert  \bigg\vert  _{\textrm{HS}} ^2\Bigg)^{1/2} =  n^{-1/2}  D_{1,n}\\
& \leq  4 \gamma_n ^{-1}n^{-1/2} \bigg(\sum_{ p \geq 1} \lambda_p ^{1/2}\bigg)\bigg(\sum_{ q \geq 1} \mu_q ^{1/2}\bigg)
\end{aligned}
\end{equation*}
from which, and Markov inequality, we deduce that    $ \big  \langle  \widetilde{V}_{XY} ^{(1)} \; , \widetilde{V}_{XY} ^{(2)} \big \rangle_{\textrm{HS}}=o_P(1)$ since $\lim_{n\rightarrow +\infty}\left(\gamma_n ^{-1}n^{-1/2}\right)=0$.
 
\bigskip

\noindent Finally, considering
\begin{equation}\label{bn}
 E_n =   \frac{1}{n} \sum_{p=1}^{+\infty} \sum_{q=1}^{+\infty}  {(\lambda_p + \gamma_n)}^{-1}{(\mu_q + \gamma_n)}^{-1}     \sum_{i=1}^{n} \bigg(\alpha_{p,i,q}^2 -  \mathbb{E}\left(\alpha_{p,i,q}^2\right)  \bigg),
\end{equation}
where $\alpha_{p,i,q}  $ is defined in \eqref{alpha}, we have:

\bigskip

\begin{lemma}\label{lem4}
Assume ($\mathscr{A}_1$) to ($\mathscr{A}_5$). Then, under $\mathscr{H}_0$, we have $E_n=o_P(1)$.
\end{lemma}
\noindent\textit{Proof}.
 Since  $E_n$ is   centered it is enough to show that  $\lim_{n\rightarrow +\infty}\left(var (E_n) \right)= 0$, then Markov inequality allows to conclude.  We have
\begin{equation*}
\begin{aligned}
var \left(E_n\right)
& =  \frac{1}{n^2}\sum_{i=1}^{n}  var \Big [   \sum_{p=1}^{+\infty} \sum_{q=1}^{+\infty}  {(\lambda_p + \gamma_n)}^{-1}{(\mu_q + \gamma_n)}^{-1}    \big  \{ \alpha_{p,i,q}^2 -  E \bigg(\alpha_{p,i,q}^2\bigg)  \big \}  \Big  ]   \\
& =  \frac{1}{n^2}\sum_{i=1}^{n}  v_{n,i},  \\
\end{aligned}
\end{equation*}
where
\begin{equation}\label{vni}
 v_{n,i} =  \sum_{p_1,q_1=1}^{+\infty} \sum_{p_2,q_2=1}^{+\infty}  {(\mu_{q_1 }+ \gamma_n)}^{-1}{(\lambda_{p_1 }+ \gamma_n)}^{-1}{(\mu_{q_2 }+ \gamma_n)}^{-1}{(\lambda_{p_2 }+ \gamma_n)}^{-1} cov ( \alpha_{p_1, i, q_1}^2  , \alpha_{p_2, i, q_2}^2  ) .
\end{equation}
Further,
\begin{equation}\label{eqt1}
\begin{aligned}
\left\vert cov ( \alpha_{p_1, i, q_1}^2  , \alpha_{p_2, i, q_2}^2  )\right\vert
 &   \leq  \mathbb{E} \bigg( \alpha_{p_1, i, q_1}^2  \alpha_{p_2, i, q_2}^2 \bigg) + \mathbb{E} \left( \alpha_{p_1, i, q_1}^2\right)\,\, \mathbb{E} \left( \alpha_{p_2, i, q_2}^2 \right)  \\
& \leq 2 \; \mathbb {E}^{1/2} \left( \alpha_{p_1, i, q_1}^4\right)   \; \mathbb{E}^{1/2} \left( \alpha_{p_2, i, q_2}^4 \right)   \\
\end{aligned}
\end{equation}
and
\begin{eqnarray}\label{helpfull3}
\alpha_{p, i, q}^2 & =& \bigg( e_p (X_i) - \mathbb{E} \left(e_p (X_i)\right) \bigg)^2 \bigg( f_q (Y_i) - \mathbb{E} \left(f_q (Y_i)\right) \bigg)^2\nonumber\\ 
&  = &\bigg(e_p^2 (X_i) + \mathbb{E}^2 \left(e_p (X_i)\right) -2e_p (X_i)\mathbb{E} \left(e_p (X_i)\right)\bigg)\nonumber\\
& &\hspace{2cm}   \,\bigg( f_q^2 (Y_i) + \mathbb{E}^2 \left(f_q (Y_i)\right)-2 f_q (Y_i) \mathbb{E} \left(f_q (Y_i)\right)\bigg) \nonumber \\
&  \leq & \bigg(e_p^2 (X_i) + \mathbb{E}^2 \left(\left\vert e_p (X_i)\right\vert \right) +2\left\vert e_p (X_i)\right\vert\mathbb{E} \left(\left\vert e_p (X_i)\right\vert\right)\bigg) \nonumber\\
& &\hspace{2cm}  \,\bigg( f_q^2 (Y_i) + \mathbb{E}^2 \left(\left\vert f_q (Y_i)\right\vert\right)+2 \left\vert f_q (Y_i)\right\vert \mathbb{E} \left(\left\vert f_q (Y_i)\right\vert \right)\bigg).   
\end{eqnarray}
Since
\begin{equation}\label{maj1}
\left\vert e_p (X_i)\right\vert=\left\vert\left\langle e_p,K(X_i,\cdot)\right\rangle_{\mathcal{H}_\mathcal{X}}\right\vert
\leq\Vert K(X_i,\cdot)\Vert_{\mathcal{H}_\mathcal{X}}= K(X_i,X_i)^{1/2}\leq \Vert K\Vert_\infty^{1/2}
\end{equation}
and, similarly,  
\begin{equation}\label{maj2}
\left\vert f_q (Y_i)\right\vert\leq \Vert L\Vert_\infty^{1/2},
\end{equation} 
 it follows from \eqref{helpfull3} that $\alpha_{p, i, q}^2
 \leq 16 \;\Vert K\Vert_{\infty} \Vert L\Vert_{\infty}$. Then, from   (\ref{eqt1}) and the equality
\begin{equation}\label{helpfull4}
\mathbb{E} \left(\alpha_{p,i,q}^2\right)=var\left( e_p (X_i)\right)\,var\left( f_q (Y_i)\right)=\left\langle e_p,\Sigma_Xe_p\right\rangle_{\mathcal{H}_\mathcal{X}}\left\langle f_q,\Sigma_Yf_q\right\rangle_{\mathcal{H}_\mathcal{Y}}=\lambda_p \mu_q
\end{equation}
we get
\begin{equation}\label{eqt5}
\begin{aligned}
\left\vert cov ( \alpha_{p_1, i, q_1}^2  , \alpha_{p_2, i, q_2}^2  )\right\vert
&   \leq   32 \;\Vert K\Vert_\infty\Vert L\Vert_\infty \lambda_{p_1 }^{1/2} \lambda_{p_2 }^{1/2}\mu_{q_1 }^{1/2}\mu_{q_2 }^{1/2}.
\end{aligned}
\end{equation}
Using this later inequality together with the inequality $(a + \gamma_n)^{-1} \leq \gamma_n^{-1}$ for all positive number $a$, we deduce from \eqref{vni} that  
$$  v_{n,i}  \leq 32 \;\Vert K\Vert_\infty\Vert L\Vert_\infty \; \gamma_n ^{-4} \Big (  \sum_{ p = 1}^{+\infty} \lambda_p ^{1/2}  \Big ) ^2 \Big (  \sum_{ q = 1}^{+\infty} \mu_q ^{1/2}  \Big ) ^2.$$
Thus
$$var \left(E_n\right) \leq 32 \;\Vert K\Vert_\infty\Vert L\Vert_\infty \; \gamma_n ^{-4} n^{-1} \Big (  \sum_{ p = 1}^{+\infty} \lambda_p ^{1/2}  \Big ) ^2 \Big (  \sum_{ q = 1}^{+\infty} \mu_q ^{1/2}  \Big ) ^2, $$
Assumptions ($\mathscr{A}_2$), ($\mathscr{A}_3$) and ($\mathscr{A}_5$) allow to conclude that $E_n = o_p(1)$.
 
\subsection{Proof of Theorem \ref{loi}}
\noindent Lemma \ref{lem2} and  an application of  Slutsky's theorem show that it is enough to prove that 
\begin{equation}\label{conv}
\widetilde{T}_n\stackrel{\mathscr{D}}{ \longrightarrow} \mathcal{N} (0,1)\,\,\,\textrm{  as }\,\,\,n\rightarrow +\infty.
\end{equation} 
From \eqref{Vchap} we have  
\[
\widetilde{T}_n=\frac{n \Vert  \widetilde{V}_{XY} ^{(1)} - \widetilde{V}_{XY} ^{(2)}  \Vert_{\textrm{HS}} ^2-D_{1,n}}{\sqrt{2}D_{2,n}}=\frac{n  \Vert  \widetilde{V}_{XY} ^{(1)}\Vert_{\textrm{HS}} ^2 -   D_{1,n}}{\sqrt{2}D_{2,n}} 
+\frac{n  \ \Vert \widetilde{V}_{XY} ^{(2)}  \Vert_{\textrm{HS}} ^2 - 2n\, \big \langle  \widetilde{V}_{XY} ^{(1)} \; , \widetilde{V}_{XY} ^{(2)} \big \rangle_{\textrm{HS}}}{\sqrt{2}D_{2,n}} .
\]
Furhermore, using Lemma \ref{lem3}-$(i)$,  (\ref{helpfull1}) and the third line  in (\ref{helpfull2}), we  obtain 
\[
n  \Vert  \widetilde{V}_{XY} ^{(1)}\Vert_{\textrm{HS}} ^2 -   D_{1,n}=n  \Vert  \widetilde{V}_{XY} ^{(1)}\Vert_{\textrm{HS}} ^2 - \mathbb{E}\bigg( n  \Vert  \widetilde{V}_{XY} ^{(1)}\Vert_{\textrm{HS}} ^2\bigg)=  E_n + 2F_n
\]
where $E_n$ is given in \eqref{bn} and
$$   F_n =  \frac{1}{n} \sum_{p=1}^{+\infty} \sum_{q=1}^{+\infty}  {(\lambda_p + \gamma_n)}^{-1}{(\mu_q + \gamma_n)}^{-1}   \Big  \{  \sum_{i=2}^{n} \alpha_{p,i,q}M_{p,i-1,q}     \Big \}  ,$$
the random variables $\alpha_{p,i,q}$ and $M_{p,i-1,q}$ being defined in \eqref{alpha} and \eqref{mpjq}.    Therefore, we can write
\[
\widetilde{T}_n=\sqrt{2}\frac{C_n}{D_{2,n}}+\frac{B_n}{\sqrt{2}D_{2,n}} 
+\frac{n  \ \Vert \widetilde{V}_{XY} ^{(2)}  \Vert_{\textrm{HS}} ^2 - 2n\, \big \langle  \widetilde{V}_{XY} ^{(1)} \; , \widetilde{V}_{XY} ^{(2)} \big \rangle_{\textrm{HS}}}{\sqrt{2}D_{2,n}} 
\]
and since,  from Lemma 18 in  \cite{harchaoui08},  we have  $\lim_{n\rightarrow +\infty}\left(D_{2,n}\right)=+\infty$, we deduce from Lemma \ref{lem3} and Lemma \ref{lem4} that
\[
\widetilde{T}_n=\sqrt{2}D_{2,n}^{-1}F_n
+o_P(1)
\]
So, in order to get \eqref{conv}  it suffices to prove that
\begin{equation}\label{ToShow}
\begin{aligned}
 D_{2,n} ^{-1} F_n \stackrel{\mathscr{D}}{ \longrightarrow} \mathcal{N} (0,1/2)\,\,\,\textrm{  as }\,\,\,n\rightarrow +\infty.
\end{aligned}
\end{equation}
We will show it   by using the central limit theorem for triangular arrays of martingale differences.
For all $(p,q)$ we set $\alpha_{p,0,q} = 0$, $ \mathscr{F}_{n,0} = \{\emptyset,\Omega  \}$ and for $i \geq 1$
$$\mathscr{F}_{n,i} = \sigma (\alpha_{p,l,q}, 1 \leq l \leq i, 1 \leq p \leq n,1 \leq q \leq n).  $$
Defining $\epsilon_{n,1}=0$ and for all $i \geq 2$,   
$$  \epsilon_{n,i} =  D_{2,n} ^{-1} n^{-1} \sum_{p=1}^{+\infty} \sum_{q=1}^{+\infty}  {(\lambda_p + \gamma_n)}^{-1}{(\mu_q + \gamma_n)}^{-1}    \alpha_{p,i,q} M_{p,i-1,q},   $$
we have  $ D_{2,n} ^{-1} F_n = \sum_{i=1}^{n} \epsilon_{n,i}$ and $\epsilon_{n,i}$ is a martingale increment since $\mathbb{E}\left(\epsilon_{n,i}\vert \mathscr{F}_{n,i-1}\right)=0$. Then, from Theorem 3.2 and Corollary 3.1 in \cite{hall80}, we will obtain \eqref{ToShow} if we show that  
\begin{equation}\label{ToVerify1}
\begin{aligned}
S_n ^2 := \sum_{i=1}^{n} \mathbb{E} \bigg( \epsilon_{n,i} ^2 |\mathscr{F}_{n,i-1} \bigg) \stackrel{ P}{ \longrightarrow}\frac{1}{2}\,\,\,\textrm{ as }\,\,\,n\rightarrow +\infty
\end{aligned}
\end{equation}
and
\begin{equation}\label{ToVerify2}
\begin{aligned}
\mathbb{E} \bigg( \max_{1 \leq i \leq n}\left( \epsilon_{n,i} ^2 \right)\bigg) = o(1).
\end{aligned}
\end{equation}

\noindent Proof of \eqref{ToVerify1}:   We have
\begin{eqnarray*}
S_n ^2 &= & D_{2,n} ^{-2} n^{-2}\sum_{i=1}^{n}  \mathbb{E} \Bigg( \bigg( \sum_{p=1}^{+\infty} \sum_{q=1}^{+\infty}  {(\lambda_p + \gamma_n)}^{-1}{(\mu_q + \gamma_n)}^{-1}    \alpha_{p,i,q} M_{p,i-1,q}\bigg)^2\bigg\vert\mathscr{F}_{n,i-1} \Bigg) \\
&=&G_n+H_n,
\end{eqnarray*}
where
\begin{eqnarray*}
 G_n =  D_{2,n} ^{-2} n^{-2} \sum_{p=1}^{+\infty} \sum_{q=1}^{+\infty} {(\mu_q + \gamma_n)}^{-2}{(\lambda_p + \gamma_n)}^{-2} \sum_{i=1}^{n} M_{p,i-1,q} ^2  \mathbb{E} \bigg(\alpha_{p,i,q}^2\bigg)
\end{eqnarray*}
and
\begin{eqnarray*}
 H_n &=&  D_{2,n} ^{-2} n^{-2}  \sum_{p_1, q_1=1}^{+\infty} \sum_{\stackrel{p_2,q_2=1}{(p_2,q_2)\neq (p_1,q_1)}}^{+\infty}\sum_{i=1}^{n}  {(\mu_{q_1 } + \gamma_n)}^{-1}{(\lambda_{p_1 } + \gamma_n)}^{-1}  \\
&&\hspace{3cm} {(\mu_{q_2 } + \gamma_n)}^{-1}{(\lambda_{p_2 } + \gamma_n)}^{-1}    M_{p_1,i-1,q_1} M_{p_2,i-1,q_2}\mathbb{E}\bigg(\alpha_{p_1, i, q_1} \alpha_{p_2, i, q_2}\bigg) .
\end{eqnarray*}
First, notice that under $\mathscr{H}_0$ one has:
\begin{eqnarray}\label{delta}
\mathbb{E}\bigg(\alpha_{p_1, i, q_1} \alpha_{p_2, i, q_2}\bigg) &=&\mathbb{E}\Bigg(\bigg(e_{p_1}(X_i)-\mathbb{E}\left(e_{p_1}(X_i)\right)\bigg)\bigg(e_{p_2}(X_i)-\mathbb{E}\left(e_{p_2}(X_i)\right)\bigg)\Bigg)\nonumber\\
&&\hspace{2cm}\,\,\mathbb{E}\Bigg(\bigg(f_{q_1}(Y_i)-\mathbb{E}\left(f_{q_1}(Y_i)\right)\bigg)\bigg(f_{q_2}(Y_i)-\mathbb{E}\left(f_{q_2}(Y_i)\right)\bigg)\Bigg)\nonumber\\
&=&\mathbb{E}\Bigg(\left\langle e_{p_1},K(X_i,\cdot)-m_X\right\rangle_{\mathcal{H}_{\mathcal{X}}}\left\langle e_{p_2},K(X_i,\cdot)-m_X\right\rangle_{\mathcal{H}_{\mathcal{X}}}\Bigg)\nonumber\\
&&\hspace{2cm}\,\,\mathbb{E}\Bigg(\left\langle f_{q_1},L(Y_i,\cdot)-m_Y\right\rangle_{\mathcal{H}_{\mathcal{Y}}}\left\langle f_{q_2},L(Y_i,\cdot)-m_Y\right\rangle_{\mathcal{H}_{\mathcal{Y}}}\Bigg)\nonumber\\
&=&\left\langle e_{p_1},\Sigma_Xe_{p_2}\right\rangle_{\mathcal{H}_{\mathcal{X}}}\left\langle f_{q_1},\Sigma_Yf_{q_2}\right\rangle_{\mathcal{H}_{\mathcal{Y}}}\nonumber\\
&=&\lambda_{p_1}\mu_{q_1}\delta_{p_1p_2}\,\delta_{q_1q_2},
\end{eqnarray}
where $\delta$ is the Kronecker symbol. Thus, $H_n=0$ and it remains to prove that $G_n$ converges in probability to $1/2$ as $n\rightarrow +\infty$. For doing that, we will show that $\lim_{n\rightarrow +\infty}\left(\mathbb{E}(G_n)\right)=1/2$ and  $G_n - \mathbb{E}(G_n) = o_p(1)$. Using the independence of the family $(\alpha_{p,i,q})_{1 \leq i \leq n}$ and the fact that these random variables are centered, we can write
\begin{equation*}
\begin{aligned}
\mathbb{E} \left( G_n \right) 
& =    D_{2,n} ^{-2} n^{-2} \sum_{p=1}^{+\infty} \sum_{q=1}^{+\infty} {(\mu_q + \gamma_n)}^{-2}{(\lambda_p + \gamma_n)}^{-2} \sum_{i=1}^{n} \mathbb{E} \left( M_{p,i-1,q} ^2\right)  \mathbb{E} \left(\alpha_{p,i,q}^2\right)\\
& =    D_{2,n} ^{-2} n^{-2} \sum_{p=1}^{+\infty} \sum_{q=1}^{+\infty} {(\mu_q + \gamma_n)}^{-2}{(\lambda_p + \gamma_n)}^{-2} \sum_{i=1}^{n} var \bigg( \sum_{\ell=0}^{i-1} \alpha_{p,\ell,q} \bigg)\,  \mathbb{E} \left(\alpha_{p,i,q}^2\right)\\
& =  D_{2,n} ^{-2} n^{-2} \sum_{p=1}^{+\infty} \sum_{q=1}^{+\infty} {(\mu_q + \gamma_n)}^{-2}{(\lambda_p + \gamma_n)}^{-2} \sum_{i=1}^{n} \sum_{\ell=0}^{i-1} \mathbb{E} \left(\alpha_{p,\ell,q} ^2 \right)  \mathbb{E} \left(\alpha_{p,i,q}^2\right).\\
\end{aligned}
\end{equation*}
using the formula
\begin{equation}\label{helpfull6}
\begin{aligned}
\bigg(\sum_{i=1}^{n} a_i\bigg)^2 = \sum_{i=1}^{n} a_i ^2 + 2  \sum_{i=1}^{n} \sum_{\ell=0}^{i-1} a_i a_\ell 
\end{aligned}
\end{equation}
and equation (\ref{helpfull4}) we obtain the equality
\begin{equation*}
\begin{aligned}
\mathbb{E} \left( G_n \right) 
& =   2 ^{-1}  D_{2,n} ^{-2} n^{-2} \sum_{p=1}^{+\infty} \sum_{q=1}^{+\infty} {(\mu_q + \gamma_n)}^{-2}{(\lambda_p + \gamma_n)}^{-2} \Bigg [   \Bigg (   \sum_{i=1}^{n}  \mathbb{E} \bigg(\alpha_{p,i,q}^2\bigg)    \Bigg  )^2     -     \sum_{i=1}^{n}  \mathbb{E}^2 \bigg(\alpha_{p,i,q}^2\bigg)           \Bigg  ] \\
& = 2 ^{-1}  D_{2,n} ^{-2} n^{-2} \sum_{p=1}^{+\infty} \sum_{q=1}^{+\infty} {(\mu_q + \gamma_n)}^{-2}{(\lambda_p + \gamma_n)}^{-2} \Bigg [   \Big (   \sum_{i=1}^{n}  \mu_q \lambda_p    \Big  )^2     -     \sum_{i=1}^{n} \mu_q ^2 \lambda_p ^2       \Bigg  ] \\
& = 2 ^{-1}  D_{2,n} ^{-2} n^{-2} \sum_{p=1}^{+\infty} \sum_{q=1}^{+\infty} {(\mu_q + \gamma_n)}^{-2}{(\lambda_p + \gamma_n)}^{-2} \mu_q ^2 \lambda_p ^2   (n^2     -  n  ) \\
& = \frac{1}{2} - \frac{1}{2n}\\
\end{aligned}
\end{equation*}
from which we deduce  that $\lim_{n\rightarrow +\infty}\left(\mathbb{E}(G_n)\right)=1/2$. For proving that $G_n - \mathbb{E}(G_n) = o_p(1)$ it is enough to show that $\lim_{n\rightarrow +\infty}\left(var(G_n)\right)=0$. We have
$$  G_n -  \mathbb  E (G_n) =  D_{2,n} ^{-2} n^{-2} \sum_{p=1}^{+\infty} \sum_{q=1}^{+\infty} {(\mu_q + \gamma_n)}^{-2}{(\lambda_p + \gamma_n)}^{-2}  Q _{n,p,q}$$
where $Q _{n,p,q} =   \sum_{i=1}^{n} \mathbb{E} \left(\alpha_{p,i,q}^2\right) \bigg(  M_{p,i-1,q} ^2 - 
 \mathbb{E} \left(M_{p,i-1,q} ^2\right)  \bigg) $.
So that
\begin{equation}\label{var}
var \left( G_n \right) =D_{2,n} ^{-4} n^{-4}\mathbb{E}\Bigg(\bigg(\sum_{p=1}^{+\infty} \sum_{q=1}^{+\infty} {(\mu_q + \gamma_n)}^{-2}{(\lambda_p + \gamma_n)}^{-2}  Q _{n,p,q}\bigg)^2\Bigg)=s_{1,n}+s_{2,n},
\end{equation}
where
\begin{equation}\label{var1}
s_{1,n} =    D_{2,n} ^{-4} n^{-4} \sum_{p=1}^{+\infty} \sum_{q=1}^{+\infty} {(\mu_q + \gamma_n)}^{-4}{(\lambda_p + \gamma_n)}^{-4} \mathbb{E} \left(Q _{n,p,q} ^2\right)
\end{equation}
and
\begin{eqnarray}\label{var2}
s_{2,n}&=  &D_{2,n} ^{-4} n^{-4}  \sum_{p_1, q_1=1}^{+\infty} \sum_{\stackrel{p_2,q_2=1}{(p_2,q_2)\neq (p_1,q_1)}}^{+\infty} {(\mu_{q_1 } + \gamma_n)}^{-2}{(\lambda_{p_1 } + \gamma_n)}^{-2}   {(\mu_{q_2 } + \gamma_n)}^{-2}{(\lambda_{p_2 } + \gamma_n)}^{-2} \nonumber\\
&&\hspace{6cm}\times\mathbb{E}\left(Q_{n,p_1,q_1}Q_{n,p_2,q_2}\right).
\end{eqnarray}
Putting
\begin{equation}\label{vijk}
\begin{aligned}
v_{p,i,q} & = M_{p,i,q} ^2 -  \mathbb{E} \left( M_{p,i,q} ^2\right) - (M_{p,i-1,q} ^2 -  \mathbb{E} \left( M_{p,i-1,q} ^2\right)) \\
&= \alpha_{p,i,q}^2- \mathbb{E} \bigg(\alpha_{p,i,q}^2\bigg) + 2 \alpha_{p,i,q} M_{p,i-1,q}.
\end{aligned}
\end{equation}
and using the formula
$$\sum_{i=1}^{n} a_i b_i= a_n \sum_{i=1}^{n}b_i - \sum_{i=1}^{n-1}\sum_{j=1}^{i}b_j (a_{i+1} - a_i)$$
with $a_i =M_{p,i-1,q} ^2 -  \mathbb{E} \left( M_{p,i-1,q} ^2\right)$ and $b_i = \mathbb{E} \bigg(\alpha_{p,i,q}^2\bigg)$   we get
\begin{equation*}
\begin{aligned}
Q _{n,p,q} &= \bigg( M_{p,n-1,q} ^2 -  \mathbb{E} \left( M_{p,n-1,q} ^2\right) \bigg)   \sum_{i=1}^{n} \mathbb{E} \bigg(\alpha_{p,i,q}^2\bigg)  -   \sum_{i=1}^{n-1} v_{p,i,q} \sum_{j=1}^{i} \mathbb{E} \left(\alpha_{p,j,q}^2\right) \\
& =  \bigg( M_{p,n-1,q} ^2 -  \mathbb{E} \left( M_{p,n-1,q} ^2\right)-  \sum_{i=1}^{n-1} v_{p,i,q} \bigg)   \sum_{i=1}^{n} \mathbb{E} \bigg(\alpha_{p,i,q}^2\bigg) \\
&+ \sum_{i=1}^{n-1} v_{p,i,q} \sum_{j=1}^{n} \mathbb{E} \bigg(\alpha_{p,j,q}^2\bigg)-   \sum_{i=1}^{n-1} v_{p,i,q} \sum_{j=1}^{i} \mathbb{E} \bigg(\alpha_{p,j,q}^2\bigg) \\
& =  \bigg( M_{p,n-1,q} ^2 -  \mathbb{E} \left( M_{p,n-1,q} ^2\right)-  \sum_{i=1}^{n-1} v_{p,i,q} \bigg) \sum_{i=1}^{n} \mathbb{E} \bigg(\alpha_{p,i,q}^2\bigg)  +   \sum_{i=1}^{n-1} v_{p,i,q} \sum_{j=i+1}^{n} \mathbb{E} \bigg(\alpha_{p,j,q}^2\bigg) \\
& = \sum_{i=1}^{n-1} v_{p,i,q} \sum_{j=i+1}^{n}\mathbb{E} \bigg(\alpha_{p,j,q}^2\bigg), \\
\end{aligned}
\end{equation*}
so that
\begin{equation*}\label{eq} 
 \mathbb{E}\bigg(Q_{n,p_1,q_1}Q_{n,p_2,q_2}\bigg) = \sum_{i=1}^{n-1}\sum_{k=1}^{n-1}\Bigg (\sum_{j=i+1}^{n} \mathbb{E} \left(\alpha_{p_1, j, q_1}^2\right) \Bigg )\Bigg ( \sum_{\ell=k+1}^{n} \mathbb{E} \left(\alpha_{p_2,\ell ,q_2}^2\right) \Bigg )  \mathbb{E} \left(v_{p_1, i ,q_1}  v_{p_2, k, q_2} \right).
\end{equation*}
From (\ref{vijk}), we get  for $i \neq k$  with $i >k$ (without loss of generality): 
\begin{equation}\label{helpfull5}
\begin{aligned}
 \mathbb{E} \left(v_{p_1, i ,q_1}  v_{p_2, k, q_2} \right) & = \mathbb{E} \Bigg( \bigg(\alpha_{p_1,i,q_1}^2- \mathbb{E} \left(\alpha_{p_1,i,q_1}^2\right) + 2 \alpha_{p_1,i,q_1} M_{p_1,i-1,q_1}\bigg)\\
&\hspace{2cm} \,\bigg(\alpha_{p_2,k,q_2}^2- \mathbb{E} \left(\alpha_{p_2,k,q_2}^2\right) + 2 \alpha_{p_2,k,q_2} M_{p_2,k-1,q_2}\bigg)\Bigg)\\
& = \mathbb{E} \Bigg( \bigg(\alpha_{p_1,i,q_1}^2- \mathbb{E} \left(\alpha_{p_1,i,q_1}^2\right) \bigg)\,\bigg(\alpha_{p_2,k,q_2}^2- \mathbb{E} \left(\alpha_{p_2,k,q_2}^2\right) \bigg)\\
& +2 \alpha_{p_1,i,q_1}\bigg(\alpha_{p_2,k,q_2}^2- \mathbb{E} \left(\alpha_{p_2,k,q_2}^2\right) \bigg)  M_{p_1,i-1,q_1}\\
&  + 2 \alpha_{p_2,k,q_2} \bigg(\alpha_{p_1,i,q_1}^2- \mathbb{E} \left(\alpha_{p_1,i,q_1}^2\right) \bigg)M_{p_2,k-1,q_2}\\
&+ 4 \alpha_{p_1,i,q_1}\alpha_{p_2,k,q_2} M_{p_1,i-1,q_1}M_{p_2,k-1,q_2} \Bigg) 
\end{aligned}
\end{equation}
and since, from independence properties, we have 
\begin{eqnarray*}
& &\mathbb{E} \Bigg( \bigg(\alpha_{p_1,i,q_1}^2- \mathbb{E} \left(\alpha_{p_1,i,q_1}^2\right) \bigg)\,\bigg(\alpha_{p_2,k,q_2}^2- \mathbb{E} \left(\alpha_{p_2,k,q_2}^2\right) \bigg)\Bigg)\\
&=&\mathbb{E} \Bigg( \bigg(\alpha_{p_1,i,q_1}^2- \mathbb{E} \left(\alpha_{p_1,i,q_1}^2\right) \Bigg)\ \,\,\mathbb{E} \Bigg( \bigg(\alpha_{p_2,k,q_2}^2- \mathbb{E} \left(\alpha_{p_2,k,q_2}^2\right) \Bigg)=0,
\end{eqnarray*}
\begin{eqnarray*}
& &\mathbb{E} \Bigg(  \alpha_{p_1,i,q_1}\bigg(\alpha_{p_2,k,q_2}^2- \mathbb{E} \left(\alpha_{p_2,k,q_2}^2\right)\bigg) M_{p_1,i-1,q_1} \Bigg) \\
&=&\mathbb{E} \bigg(  \alpha_{p_1,i,q_1}\bigg)\,\,\mathbb{E} \Bigg(  \bigg(\alpha_{p_2,k,q_2}^2- \mathbb{E} \left(\alpha_{p_2,k,q_2}^2\right)\bigg) M_{p_1,i-1,q_1} \Bigg)=0,
\end{eqnarray*}
\begin{eqnarray*}
&&\mathbb{E} \Bigg(  \alpha_{p_2,k,q_2} \bigg(\alpha_{p_1,i,q_1}^2- \mathbb{E} \left(\alpha_{p_1,i,q_1}^2\right) \bigg)M_{p_2,k-1,q_2} \Bigg) \\
&=&\mathbb{E} \bigg(  \alpha_{p_2,k,q_2} \bigg)\mathbb{E} \Bigg( \alpha_{p_1,i,q_1}^2- \mathbb{E} \left(\alpha_{p_1,i,q_1}^2\right) \bigg)M_{p_2,k-1,q_2} \Bigg) =0,
\end{eqnarray*}
\[
 \mathbb{E} \Bigg(  \alpha_{p_1,i,q_1}\alpha_{p_2,k,q_2} M_{p_1,i-1,q_1}M_{p_2,k-1,q_2} \Bigg) = \mathbb{E} \bigg(  \alpha_{p_1,i,q_1}\bigg) \mathbb{E} \Bigg( \alpha_{p_2,k,q_2} M_{p_1,i-1,q_1}M_{p_2,k-1,q_2} \Bigg)=0,
\]
it follows that  $ \mathbb{E} \left(v_{p_1, i ,q_1}  v_{p_2, k, q_2} \right)=0$. Thus
\begin{equation*}\label{eq2} 
 \mathbb{E}\bigg(Q_{n,p_1,q_1}Q_{n,p_2,q_2}\bigg) = \sum_{i=1}^{n-1}\Bigg (\sum_{j=i+1}^{n} \mathbb{E} \left(\alpha_{p_1, j, q_1}^2\right) \Bigg )\Bigg ( \sum_{\ell=i+1}^{n} \mathbb{E} \left(\alpha_{p_2, \ell ,q_2}^2\right) \Bigg )  \mathbb{E} \left(v_{p_1, i ,q_1}  v_{p_2, i, q_2} \right)
\end{equation*}
and
\begin{eqnarray}\label{eq3} 
 \bigg\vert\mathbb{E}\bigg(Q_{n,p_1,q_1}Q_{n,p_2,q_2}\bigg) \bigg\vert&\leq& \Bigg (\sum_{j=1}^{n} \mathbb{E} \left(\alpha_{p_1, j, q_1}^2\right) \Bigg )\Bigg ( \sum_{\ell=1}^{n} \mathbb{E} \left(\alpha_{p_2, \ell ,q_2}^2\right) \Bigg )  \sum_{i=1}^{n-1}\bigg\vert\mathbb{E} \left(v_{p_1, i ,q_1}  v_{p_2, i, q_2} \right)\bigg\vert\nonumber\\
&  \leq  &  n^2  \mu_{q_1 } \lambda_{p_1 } \mu_{q_2 } \lambda_{p_2 }  \sum_{i=1}^{n-1}\bigg\vert\mathbb{E} \left(v_{p_1, i ,q_1}  v_{p_2, i, q_2} \right)\bigg\vert .
\end{eqnarray}
As in \eqref{helpfull5} we get 
\begin{eqnarray}\label{ev}
 \mathbb{E} \left(v_{p_1, i ,q_1}  v_{p_2, i, q_2} \right) 
& = &\mathbb{E} \Bigg( \bigg(\alpha_{p_1,i,q_1}^2- \mathbb{E} \left(\alpha_{p_1,i,q_1}^2\right) \bigg)\,\bigg(\alpha_{p_2,i,q_2}^2- \mathbb{E} \left(\alpha_{p_2,i,q_2}^2\right) \bigg)\nonumber\\
&&+2 \alpha_{p_1,i,q_1}\bigg(\alpha_{p_2,i,q_2}^2- \mathbb{E} \left(\alpha_{p_2,i,q_2}^2\right) \bigg)  M_{p_1,i-1,q_1}\nonumber\\
& &  + 2 \alpha_{p_2,i,q_2} \bigg(\alpha_{p_1,i,q_1}^2- \mathbb{E} \left(\alpha_{p_1,i,q_1}^2\right) \bigg)M_{p_2,i-1,q_2})\nonumber\\
&&+ 4 \alpha_{p_1,i,q_1}\alpha_{p_2,i,q_2} M_{p_1,i-1,q_1}M_{p_2,i-1,q_2} \Bigg) \nonumber\\
& =&cov ( \alpha_{p_1, i, q_1}^2  , \alpha_{p_2, i, q_2}^2  ) +2\mathbb{E} \Bigg(  \alpha_{p_1,i,q_1}\bigg(\alpha_{p_2,i,q_2}^2- \mathbb{E} \left(\alpha_{p_2,i,q_2}^2\right) \bigg)\Bigg)\,\mathbb{E} \bigg(   M_{p_1,i-1,q_1}\bigg)\nonumber\\
&  &+ 2\mathbb{E} \Bigg(\alpha_{p_2,i,q_2} \bigg(\alpha_{p_1,i,q_1}^2- \mathbb{E} \left(\alpha_{p_1,i,q_1}^2\right) \bigg)\Bigg)\,\mathbb{E} \bigg(M_{p_2,i-1,q_2}\bigg))\nonumber\\
&&+ 4 \mathbb{E} \bigg(\alpha_{p_1,i,q_1}\alpha_{p_2,i,q_2} \bigg)\,\mathbb{E} \bigg(M_{p_1,i-1,q_1}M_{p_2,i-1,q_2} \bigg)\nonumber \\
& =& cov ( \alpha_{p_1, i, q_1}^2  , \alpha_{p_2, i, q_2}^2  ) + 4 \mathbb{E} \bigg(\alpha_{p_1,i,q_1}\alpha_{p_2,i,q_2} \bigg)\,\mathbb{E} \bigg(M_{p_1,i-1,q_1}M_{p_2,i-1,q_2} \bigg) .
\end{eqnarray}
Moreover the independence of the variables $\alpha_{p_1,i,q_1}$ and $\alpha_{p_2,j,q_2}$ for $i\neq j$ gives  
\begin{equation*}\label{helpfull9}
\begin{aligned}
\mathbb{E} \bigg(M_{p_1,i-1,q_1}M_{p_2,i-1,q_2} \bigg)=  \sum_{\ell=0}^{i-1} \mathbb{E} \bigg(\alpha_{p_1, \ell ,q_1 }\alpha_{p_2,\ell ,q_2 }\bigg)
\end{aligned}
\end{equation*}
and  using again (\ref{helpfull6}) and \eqref{delta}, we obtain
\begin{eqnarray}\label{somme}
&&\sum_{i=1}^{n-1}   \bigg\vert     \mathbb{E} \bigg(\alpha_{p_1, i ,q_1 }\alpha_{p_2,i,q_2 }\bigg)\,\mathbb{E} \bigg(M_{p_1,i-1,q_1}M_{p_2,i-1,q_2} \bigg)          \bigg\vert    \nonumber\\
& =& \sum_{i=1}^{n-1}   \bigg\vert    \mathbb{E} \bigg(\alpha_{p_1, i ,q_1 }\alpha_{p_2,i,q_2 }\bigg)\sum_{\ell=0}^{i-1} \mathbb{E} \bigg(\alpha_{p_1, \ell ,q_1 }\alpha_{p_2,\ell ,q_2 }\bigg)         \bigg\vert  \nonumber  \\
& \leq &\frac{1}{2}  \Bigg(   \Bigg (   \sum_{i=1}^{n-1}   \bigg\vert   \mathbb{E} \bigg(\alpha_{p_1, i ,q_1 }\alpha_{p_2,i,q_2 }\bigg)\bigg\vert   \,  \Bigg  )^2     -     \sum_{i=1}^{n-1}  \mathbb{E}^2 \bigg(\alpha_{p_1, i ,q_1 }\alpha_{p_2,i,q_2 }\bigg)          \Bigg  )\nonumber \\
& \leq &\frac{1}{2}    \Bigg (   \sum_{i=1}^{n-1}   \bigg\vert   \mathbb{E} \bigg(\alpha_{p_1, i ,q_1 }\alpha_{p_2,i,q_2 }\bigg)\bigg\vert   \,  \Bigg  )^2\nonumber \\
& \leq &\frac{1}{2} n^2 \mu_{q_1 } \lambda_{p_1 }\mu_{q_2 } \lambda_{p_2 } \delta_{p_1 p_2 }\delta_{q_1 q_2 }.
\end{eqnarray}
Finally, \eqref{ev}, (\ref{eqt5}) and \eqref{somme} give
$$  \sum_{j=1}^{n-1}  \bigg\vert          \mathbb{E} [ v_{p_1, j ,q_1}  v_{p_2, j ,q_2}   ]    \bigg\vert   \leq   32 n\;\Vert K\Vert_\infty\Vert L\Vert_\infty \lambda_{p_1 }^{1/2} \lambda_{p_2 }^{1/2}\mu_{q_1 }^{1/2}\mu_{q_2 }^{1/2} + 2 n^2 \mu_{q_1 } \lambda_{p_1 }\mu_{q_2 } \lambda_{p_2 }  \delta_{p_1 p_2 }\delta_{q_1 q_2 }  $$
and from \eqref{eq3} we obtain:
\begin{equation}\label{eq4}
  \bigg\vert\mathbb{E}\bigg(Q_{n,p_1,q_1}Q_{n,p_2,q_2}\bigg) \bigg\vert \leq 2n^4    \bigg( 16n^{-1}  \Vert K\Vert_\infty\Vert L\Vert_\infty \;  \mu_{q_1 } ^{3/2} \lambda_{p_1 }^{3/2}\mu_{q_2 }^{3/2} \lambda_{p_2 }^{3/2} +  \mu_{q_1 }^2 \lambda_{p_1 }^2\mu_{q_2 }^2 \lambda_{p_2 }^2 \delta_{p_1 p_2 }\delta_{q_1 q_2 }     \bigg).      
\end{equation}
This inequality allows to obtain bounds for the terms $s_{1,n}$ and $s_{2,n}$ given in \eqref{var1} and \eqref{var2}. Indeed, from \eqref{eq4}, we deduce that
\begin{equation*}\label{v1}
 \mathbb{E}\bigg(Q_{n,p,q}^2\bigg)  \leq 2n^4    \bigg( 16n^{-1}  \Vert K\Vert_\infty\Vert L\Vert_\infty \;  \lambda_{p }^{3}\mu_{q}^{3} +    \lambda_{p}^4\mu_{q }^4     \bigg) 
\end{equation*}
and, therefore, that
\begin{eqnarray*}
s_{1,n} &\leq&    2D_{2,n} ^{-4} \Bigg\{16n^{-1}  \Vert K\Vert_\infty\Vert L\Vert_\infty \sum_{p=1}^{+\infty}\frac{\lambda_p^3}{(\lambda_p + \gamma_n)^{4} } \sum_{q=1}^{+\infty}\frac{\mu^3_q}{(\mu_q + \gamma_n)^{4}}\\
&&+\sum_{p=1}^{+\infty}\left(\frac{\lambda_p}{\lambda_p + \gamma_n}\right)^4 \sum_{q=1}^{+\infty}\left(\frac{\mu_q}{\mu_q + \gamma_n}\right)^4\Bigg\}.
\end{eqnarray*}
Since
\[
\frac{\lambda_p^3}{(\lambda_p + \gamma_n)^{4} }\leq \gamma^{-1}_n\frac{\lambda_p^3}{(\lambda_p + \gamma_n)^{3} }\leq \gamma^{-1}_n\left(\frac{\lambda_p}{\lambda_p + \gamma_n}\right)^2,\,\,\left(\frac{\lambda_p}{\lambda_p + \gamma_n}\right)^4\leq \left(\frac{\lambda_p}{\lambda_p + \gamma_n}\right)^2
\]
and, similarly,
\[ 
\frac{\mu^3_q}{(\mu_q + \gamma_n)^{4}}\leq \gamma^{-1}_n\left(\frac{\mu_q}{\mu_q + \gamma_n}\right)^2\,\,\textrm{ and }\,\,\left(\frac{\mu_q}{\mu_q + \gamma_n}\right)^4\leq \left(\frac{\mu_q}{\mu_q + \gamma_n}\right)^2,
\]
it follows
\begin{eqnarray*}
s_{1,n} &\leq&    2D_{2,n} ^{-4}\sum_{p=1}^{+\infty}\left(\frac{\lambda_p}{\lambda_p + \gamma_n}\right)^2 \sum_{q=1}^{+\infty}\left(\frac{\mu_q}{\mu_q + \gamma_n}\right)^2 \Bigg\{16 \gamma^{-2}_nn^{-1}  \Vert K\Vert_\infty\Vert L\Vert_\infty 
+1\Bigg\}\\
&&= 2D_{2,n} ^{-2} \Bigg\{16 \gamma^{-2}_nn^{-1}  \Vert K\Vert_\infty\Vert L\Vert_\infty 
+1\Bigg\}
\end{eqnarray*}
and using Lemma 18 in  \cite{harchaoui08}   and Assumption ($\mathscr{A}_5$) we deduce that $\lim_{n\rightarrow +\infty}(s_{1,n})=0$. On the other hand,  \eqref{eq4} implies that  for $(p_2,q_2)\neq (p_1,q_1)$:
\begin{equation*}
  \bigg\vert\mathbb{E}\bigg(Q_{n,p_1,q_1}Q_{n,p_2,q_2}\bigg) \bigg\vert \leq 32n^3    \Vert K\Vert_\infty\Vert L\Vert_\infty \;  \mu_{q_1 } ^{3/2} \lambda_{p_1 }^{3/2}\mu_{q_2 }^{3/2} \lambda_{p_2 }^{3/2}    .      
\end{equation*}
Hence
\begin{eqnarray*}
\vert s_{2,n}\vert \leq  32D_{2,n} ^{-4} n^{-1} \Vert K\Vert_\infty\Vert L\Vert_\infty\Bigg(  \sum_{p=1}^{+\infty}\frac{\lambda_p^{3/2}}{(\lambda_p + \gamma_n)^{2} }\Bigg)^2\,\Bigg( \sum_{q=1}^{+\infty}\frac{\mu_q^{3/2}}{(\mu_q + \gamma_n)^{2} }\Bigg)^2.
\end{eqnarray*}
Since
\[
\frac{\lambda_p^{3/2}}{(\lambda_p + \gamma_n)^{2} }=\frac{\lambda_p}{\lambda_p + \gamma_n }\frac{\lambda_p^{1/2}}{\lambda_p + \gamma_n }\leq \frac{\lambda_p^{1/2}}{\lambda_p + \gamma_n }\leq \gamma^{-1}_n\lambda_p^{1/2}\,\,\textrm{ and   }\,\,\frac{\mu_q^{3/2}}{(\mu_q + \gamma_n)^{2} }\leq \gamma^{-1}_n\mu_q^{1/2},
\]
it follows
\begin{eqnarray*}
\vert s_{2,n}\vert \leq  32D_{2,n} ^{-4} \gamma^{-4}_nn^{-1} \Vert K\Vert_\infty\Vert L\Vert_\infty\Bigg(  \sum_{p=1}^{+\infty}\lambda_p^{1/2}\Bigg)^2\,\Bigg( \sum_{q=1}^{+\infty}\mu_q^{1/2}\Bigg)^2.
\end{eqnarray*}
Lemma 18 in  \cite{harchaoui08}   and Assumptions ($\mathscr{A}_2$), ($\mathscr{A}_3$) and  ($\mathscr{A}_5$) allow to  deduce that $\lim_{n\rightarrow +\infty}(s_{2,n})=0$. Finally, from \eqref{var} we conclude that $\lim_{n\rightarrow +\infty}(var (G_n)) =0$, what  ensures that  $G_n -\mathbb{E}(G_n)= o_P(1)$.

\bigskip

\noindent Proof of \eqref{ToVerify2}: Using   \eqref{maj1} et \eqref{maj2}  we obtain  $\vert\alpha_{p,i,q}\vert \leq 4\Vert K\Vert_\infty^{1/2}\Vert L\Vert_\infty^{1/2}$. Thus

$$\max_{1 \leq i \leq n} |\epsilon_{n,i}| \leq 4 D_{2,n} ^{-1} n^{-1}\Vert K\Vert_\infty^{1/2}\Vert L\Vert_\infty^{1/2} \sum_{p=1}^{+\infty} \sum_{q=1}^{+\infty}  {(\lambda_p + \gamma_n)}^{-1}{(\mu_q + \gamma_n)}^{-1} \max_{1 \leq i \leq n}   | M_{p,i-1,q}|   ,$$
and from the Minkowski inequality:
$$ \mathbb{E}^{1/2} \Big [ \max_{1 \leq i \leq n} \epsilon_{n,i} ^2\Big ] \leq  4D_{2,n} ^{-1} n^{-1}\Vert K\Vert_\infty^{1/2}\Vert L\Vert_\infty^{1/2}\sum_{p=1}^{+\infty} \sum_{q=1}^{+\infty}  {(\lambda_p + \gamma_n)}^{-1}{(\mu_q + \gamma_n)}^{-1}\mathbb{E}^{1/2} \Big [ \max_{1 \leq i \leq n}   | M_{p,i-1,q}|^2\Big  ]  .   $$
One easily verifies that $(M_{p,i,q})_{1 \leq i \leq n }$ is a $\mathcal F_{n,i}$-martingale. Doob's inequality gives
$$ \mathbb{E}^{1/2} \Big [ \max_{1 \leq j \leq n}   | M_{p,i-1,q}|^2\Big  ] \leq   \mathbb{E}^{1/2} \Big [   M_{i,n-1,k}^2\Big  ] \leq n^{1/2} \mu_q ^{1/2} \lambda_p ^{1/2};  $$
so that
$$\mathbb{E } \Big [ \max_{1 \leq i \leq n} |\epsilon_{n,i}| \Big ] \leq \mathbb{E}^{1/2} \Big [ \max_{1 \leq j \leq n} \epsilon_{n,i} ^2\Big ] \leq4  D_{2,n} ^{-1} n^{-1/2} \gamma_n^{-2}\Vert K\Vert_\infty^{1/2}\Vert L\Vert_\infty^{1/2}\Big( \sum_{p=1}^{+\infty}\lambda_p ^{1/2}\Big )\Big(\sum_{q=1}^{+\infty} \mu_q ^{1/2}\Big) . $$
Assumptions ($\mathscr{A}_2$), ($\mathscr{A}_3$) and  ($\mathscr{A}_5$)  allow to conclude (\ref{ToVerify2}).



\medskip




\section*{}

\end{document}